\documentclass[a4paper, 11pt]{article}
\usepackage{fullpage}
\usepackage{geometry}
\usepackage{amsmath, amssymb}
\usepackage{amssymb}
\usepackage{epsfig}
\usepackage{amsfonts}
\usepackage{amsmath}
\usepackage{euscript}
\usepackage{amscd}
\usepackage{amsthm}
\usepackage{tikz-cd}
\usepackage{theoremref}
\usepackage{enumerate}
\usepackage{array}
\usepackage{mathrsfs}
\usepackage{tikz-cd}
\usepackage{adjustbox}
\usepackage{mathtools}
\usepackage{url}

\DeclareMathAlphabet{\mathpzc}{OT1}{pzc}{m}{it}
\DeclareMathAlphabet{\mathpzc}{OT1}{pzc}{m}{it}
\newtheorem{thm}{Theorem}[section]
\newtheorem{lem}[thm]{Lemma}
\newtheorem{prop}[thm]{Proposition} 
\newtheorem{cor}[thm]{Corollary}
\newtheorem{rem}[thm]{Remark}

\newcommand{\bZ}{\mathbb Z}

\newcommand{\A}{\mathbb A}

\newcommand{\Spec}{\operatorname{Spec}}

\newcommand{\ml}{\operatorname{ML}}

\title{On $\A^1$-contractibility of certain simple birational extensions of affine spaces} 
 \author{
 Parnashree Ghosh\\
\small{\it Department of Mathematics,}\\
\small{\it Indian Institute of Science,}\\
\small{\it CV Raman Avenue, Bengaluru-560012, India.}\\
\small{\it email: parnashreeg@iisc.ac.in, ghoshparnashree@gmail.com}
 }

\begin{document}
	
	\date{}
		\maketitle
	
		\begin{abstract}
        Over a field of characteristic zero, upto isomorphism of varieties, affine spaces are the only smooth $\A^1$-contractible affine varieties in dimensions $\leqslant 2$. However, in dimensions $\geqslant 3$, examples of smooth $\A^1$-contractible affine varieties, not isomorphic to affine spaces are constructed in \cite{DF} and \cite{DG}. 
In this paper, we consider a generalized class of smooth affine varieties containing the examples of $\A^1$-contractible varieties in \cite{DF} and \cite{DG}, given by
$$
a(x_m)b(x_1,\ldots,x_{m-1})y+f(z,t)+x_m=0,
$$
and investigate when such a variety is actually isomorphic to an affine space. We establish that, for a large subfamilies of these varieties to be affine spaces, it is necessary and sufficient that the embedding of the corresponding hyperplane in $\A^{m+3}$ must be rectifiable, in the sense that, there exists an automorphism of the ambient space that takes it to a coordinate hyperplane. Thus our result naturally connects $\mathbb A^1$-contractibility of these varieties with the classical embedding problem for affine spaces in codimension one, and provide new evidences towards the conjecture of Abhyankar and Sathaye. A key ingredient in our approach is to study singular $\mathbb A^1$-contractible affine curves over perfect fields. We describe the possible singularities of such curves and obtain a generalization of the classical result of Lin–Zaidenberg for topologically contractible affine plane curves. We further construct new examples of $\mathbb{A}^1$-contractible affine varieties in every dimension $\geqslant 4$, in the sense that they are not isomorphic to the existing examples.

	\end{abstract}



	\section{Introduction}
    Throughout this paper, unless specified, $k$ denotes a field and $\overline{k}$ denotes its algebraic closure. For a ring $R$, $R^{[n]}$ denotes a polynomial ring in $n$ variables over $R$ and $R^*$ denotes the group of units in $R$.

   Affine spaces $\A_k^n$ are the most fundamental objects in Algebraic geometry. One of the major open problems related to affine spaces is to characterize them among smooth $n$-dimensional affine varieties.  
   When $k \subset \mathbb{C}$, the affine spaces satisfies a very nice property that their underlying complex analytic space is topologically contractible. Therefore it is natural to ask whether topological contractibility is enough to characterize affine spaces among smooth affine $\mathbb{C}$-varieties. Answer to this question is affirmative in dimension 1, however in higher dimensions there are examples of topologically contractible smooth affine varieties which are not isomorphic to affine spaces. In dimension two, the Ramanujam surface (\cite{Ra}) and tom Dieck Petrie surfaces (\cite{tdp}) are examples of such smooth affine surfaces. In dimensions $\geqslant 3$, examples of topologically contractible smooth affine $\mathbb{C}$-varieties which are not affine spaces are given by Koras-Russell threefolds of first kind and their higher dimensional analogues given by
   the equations 
   \begin{equation}\label{k}
   x_1^{r_1}\cdots x_m^{r_m}y+z^r+t^s+x_m=0
   \end{equation}
   in $\A^{m+3}$ with $r_i,r,s \geqslant 2$ and $\mathrm{gcd}(r,s)=1$. The fact that these varieties are not isomorphic to affine spaces is established by studying their Makar-Limanov invariants related to the $\mathbb{G}_a$ actions on the corresponding varieties (\cite{ML}, \cite{adv1}, \cite{adv2}).
   These varieties are extremely significant in the context of Zariski Cancellation Problem in characteristic zero, which asks whether an $n$-dimensional smooth affine variety $X$ satisfying $X \times \A^1 \cong \A^{n+1}$ is isomrphic to $\A^n$. 
   Being topologically contractible, the underlying topological manifolds of these varieties are diffeomorphic to $\mathbb{R}^{2m+4}$, but it is still not known whether the cylinders over these varieties are isomorphic to affine spaces. This makes them potential counter examples to this problem. Their candidature as counter examples to this problem is further strengthened, as they turned out to be contractible in the sense of $\A^1$-homotopy theory (\cite{DF}, \cite{DG}), introduced by Morel and Voevodsky in \cite{MV}. The notion of contractibility in $\A^1$-homotopy category emerged as a stronger notion of contractibility of varieties than topological contractibility, as firstly, this gives a notion of contractibility of varieties over any field, and secondly, if $k$ has an embedding in $\mathbb{C}$, $\A^1$-contractibility of a $k$-variety implies that it must be topologically contractible. We now briefly describe the notion of $\A^1$-contractibility: 
    Let $Sm_k$ be the category of smooth, separated, finite type $k$-schemes. A $k$-scheme $X$ induces a Nisnevich sheaf $h_X$ on $Sm_k$, whose sections are given by $h_X(U):=Hom(U,X)$, for $U \in Sm_k$. Let $\mathcal{H}((Sm_k)_{Nis},\A^1)$ denotes the  unstable $\A^1$-homotopy category associated to $\Delta^{op}Shv_{Nis}(Sm_k)$.
    A $k$-scheme $X$ is said to be $\A^1$-contractible in $\mathcal{H}((Sm_k)_{Nis},\A^1)$ if the structure morphism $X \to \mathrm{Spec}(k)$ induces an isomorphism in $\mathcal{H}((Sm_k)_{Nis},\A^1)$, and $X$ is said to be stably $\A^1$-contractible if it becomes $\A^1$-contractible after $\mathbb{P}^1$-stabilization. 
    The affine spaces $\A_k^n$ are the most elementary examples of $\A^1$-contractible algebraic varieties. 
    Thus it is natural to investigate whether $\A^1$-contractibility is sufficient to characterize the affine spaces among smooth affine $k$-varieties. The answer is affirmative for smooth affine varieties upto dimension two over fields of characteristic zero(\cite{CR}). Whereas in higher dimensions the varieties given by  
    $$
    \mathrm{Spec}\,B_m:=\dfrac{k[x_1,\ldots,x_m,y,z,t]}{(x_1^{r_1}\cdots x_m^{r_m}y+z^r+t^s+x_m)}, \text{~where~}r_i,r,s \geqslant 2, \text{~with~} \mathrm{gcd}(r,s)=1,
    $$
    are proven to be $\A^1$-contractible smooth affine varieties not isomorphic to affine spaces (\cite{DF}, \cite{DG}) over any field of characteristic zero. This motivates us to ask the following question for a larger family of varieties which contains $\mathrm{Spec}\,B_m$:

     \medskip
    \noindent
    {\bf Question 1.}
    Let 
    \begin{equation}\label{a}
A:=
\dfrac{k[x_1,\ldots,x_m,y,z,t]}
{\left( a(x_m)b(x_1,\ldots,x_{m-1})y+f(z,t)+x_m\right)}, \quad a,b \neq 0, \text{ and } f \notin k
\end{equation}
defines a smooth affine variety. 

\begin{enumerate}[\rm (i)]
    \item When is $\mathrm{Spec}\,A$ $\A^1$-contractible?

    \item When is $A=k^{[m+2]}$, i.e., $\mathrm{Spec}\,A \cong \A_k^{m+2}$?
\end{enumerate}

   Before discussing on Question 1(i), we want to emphasize on the conditions for topological contractibility of $\mathrm{Spec}\,A$ over $\mathbb{C}$, which is established in the work of Kaliman and Zaidenberg in \cite{afmod} on simple birational extensions of topologically contractible $\mathbb{C}$-varieties.  We note that the variety $\mathrm{Spec}\,A$ (as in \eqref{a}) is a simple birational extension of $X:=\mathrm{Spec}\,(k[x_1,\ldots,x_m,z,t])$ with respect to center $\mathcal{Z}=V\left(a(x_m)b(x_2,\ldots,x_{m-1}), f(z,t)+x_m)\right)$ and divisor $\mathcal{D}=V\left(a(x_m)b(x_1,\ldots,x_{m-1})\right)$, that means $$\mathrm{Spec}\,A=\mathrm{Bl}_{\mathcal{Z}}X \setminus \mathcal{D}^{pr}.$$ 
   Let $\mathcal{E_Z}$ be the exceptional divisor of the corresponding blow up. Then $\mathcal{E}:=\mathcal{E_Z} \setminus \mathcal{D}^pr$ is caled the exceptional divisor of the extension.
   In \cite{afmod} Kaliman and Zaidenberg proved the following which describes the sufficient condition for $\mathrm{Spec}\,A$ to be topologically contractible.

 \begin{thm}(\cite[Theorem 3.1]{afmod})
     Let $X$ be a topologically contractible affine variety and $\sigma:X^{\prime} \to X$ be a simple birational extension of $X$ with respect to center $\mathcal{Z}$ and divisor $\mathcal{D}$. Let $\mathcal{E}$ be the exceptional divisor of the extension. If  $\mathcal{E}=\cup_{i=1}^r \mathcal{E}_i,$ and $\mathcal{D}=\cup_{i=1}^r\mathcal{D}_i$ be decomposition into irreducible components such that each component is a topological manifolds,  $\sigma^*(\mathcal{D}_i)=\mathcal{E}_i$ for every $i$ and $\mathcal{Z} \hookrightarrow \mathcal{D}$ is an weak equivalence then $X^{\prime}$ is topologically contractible.
 \end{thm}
 
We now turn to contractibility in the $\A^1$-homotopy category. The stable $\A^1$-contractibility of the Koras--Russell threefolds $\mathrm{Spec}\,B_1$, which belong to the family of varieties of the form $\mathrm{Spec}\,A$, relies on the fact that the higher Chow groups of the center
$$
\mathcal{Z}=\operatorname{Spec}\bigl(k[z,t]/(z^s+t^u)\bigr),
$$
of the corresponding simple birational extension of $\mathrm{Spec}(k[x_1,z,t])$, coincide with those of the base field (\cite{HKO}). This property is also crucial for proving $\A^1$-contractibility of Koras--Russell threefolds in the unstable setting (\cite{DF}). In light of these observations, we investigate the following sub-question, which addresses Question 1(i).

\smallskip
\noindent
\textbf{Question 1.1.}
Let $A$ be as in \eqref{a}, and suppose that $\mathrm{Spec}\,A$ is a smooth affine $k$-variety. Let $\mathcal{Z}$ and $\mathcal{D}$ be the center and divisor of the corresponding birational extension of $\operatorname{Spec}k[x_1,\dots,x_m,z,t]$. Is $\A^1$-weak equivalence of
$$
\mathcal{Z}\hookrightarrow \mathcal{D}
$$
necessary and sufficient for the $\A^1$-contractibility of $\mathrm{Spec}\,A$?

\smallskip
Equivalently, one may ask whether, for every root $\lambda$ of $a(x_m)$ in $\overline{k}$, $\A^1$-contractibility of the curves
$$
\Gamma_{\lambda}:=\mathrm{Spec}(k(\lambda)[z,t]/(f+\lambda))
$$
in the corresponding unstable $\A^1$-homotopy category is necessary and sufficient for the $\A^1$-contractibility of $\mathrm{Spec}\,A$.

We show that, over perfect fields, the answer to the necessary part of Question 1.1 is affirmative. More precisely, we prove the following result (\thref{thm3}).

\begin{thm}\thlabel{1.2}
   Let 	$A$ be a ring as in as in \eqref{a},   $\lambda$ be a root of $a(x_m)$ in $\overline{k}$, $\Gamma_{\lambda}:=\mathrm{Spec}\,(k(\lambda)[z,t]/(f+\lambda))$ and $\overline{\Gamma}_{\lambda}:=\Gamma_{\lambda} \times_{k(\lambda)} \overline{k}$.  
            If $\mathrm{Spec}\, A$ is a smooth $\mathbb{A}^1$-contractible affine variety. Then

            \begin{itemize}
                \item[(a)]
             $\Gamma_{\lambda}$ has only unibranched singularity and the normalization of $\overline{\Gamma}_{\lambda}$ is isomorphic to $\mathbb{A}_{\overline{k}}^1$.  Furthermore, the corresponding Nisnevich sheaves  $h_{\overline{\Gamma}_{\lambda}}, h_{\mathbb{A}_{\overline{k}}^1}$ over $Sm_{\overline{k}}$ are isomorphic. 
             
                \item[(b)] If $k$ is a perfect field, then the curve $\Gamma_{\lambda}$ is a polynomial curve, that means the normalization of $\Gamma_{\lambda}$ is isomorphic to $\mathbb{A}_{k(\lambda)}^1$. Furthermore, the corresponding Nisnevich sheaves $h_{\Gamma_{\lambda}}, h_{\mathbb{A}_{k(\lambda)}^1}$ over $Sm_{k(\lambda)}$ are isomorphic.
            \end{itemize}
        \end{thm}

        We use techniques from G-theory and higher Chow groups to prove the above result. 
        As an application we construct examples of affine varieties which are $\mathbb{A}^1$-chain connected but not $\mathbb{A}^1$-contractible. Also in \cite{GGPJA} several examples of factorial affine varieties are constructed that are not affine spaces. Here we show that they are not even $\mathbb{A}^1$-contractible.
        We further show that over algebraically closed fields of characteristic zero, $\A^1$-contractibility of the curves $\Gamma_{\lambda}$ in $\mathcal{H}((Sm_k)_{Nis}, \A^1)$ are sufficient to ensure stable $\A^1$-contractibility of $\mathrm{Spec}\, A$ (\thref{sufficient}), and for a subfamily of these varieties the criterion is indeed sufficient for $\A^1$-contractibility in the unstable setting (Theorem~\ref{sufficientunstable}).
        In order to prove this sufficiency part, we studied possible singularities of an $\A^1$-contractible affine curve and prove the following result (\thref{uni-branched},\thref{A1}). 
\begin{thm}
            Let $k$ be a perfect field. A geometrically integral, $\A^1$-contractible  affine $k$-curve $X$ has at most unibranched singularity.  Furthermore, the normalization of $X$ is $\A_k^1$ and $X$ is isomorphic to $\A_k^1$ as Nisnevich sheaf of sets on $Sm_k$.
        \end{thm}
        
        The above result is an extension of the classical result of Lin-Zaidenberg which states that a topologically contractible integral affine plane curve over $\mathbb{C}$ can have at most uni cuspidal singularity (\cite{LZ}). 

     Next we answer Question 1(ii) and show that our results connects to the study of Embeddings/Epimorphism  problem between affine spaces in codimension one, and the related conjecture of Abhyankar and Sathaye on ``rectifiability" of embeddings of $\A_k^{n-1}$ in $\A_k^n$. 
To be more precise, Embedding/Epimorphism problem asks the following:
	
	\smallskip
	\noindent
	\textbf{Question 2:} Let $k$ be a field and $k[x_1, \ldots, x_n]/(H) \cong k^{[n-1]}$. Is then
	\[
	k[x_1, \ldots, x_n] = k[H]^{[n-1]}?
	\]
	
	\smallskip
    If the answer to the above question is affirmative, then we say $H=0$ defines a rectifiable embedding of $\A_k^{n-1}$ in $\A_k^n$.
    The famous Epimorphism theorem of Abhyankar-Moh gives affirmative answer to Question 2 for $n=2$ and over fields of chracteristic zero (\cite{AM}), whereas over fields of positive characteristic, counter examples to Question 2 exists for every $n \geqslant 2$ (\cite{Se}, \cite{Na}). 
    The conjecture of Abhyankar and Sathaye asserts an affirmative answer to Question 2 when $n>2$ and $k$ is a field of characteristic zero. Now answering Question 1(ii) we prove the following (Theorems \ref{embedding1} and \ref{embedding}): 
    
  \begin{thm}\thlabel{embedding1}
        Let $k$ be a field of characteristic zero, and 
        $$
        A=\dfrac{k[x_1,\ldots,x_m,y,z,t]}{(H)},
        $$
        where
        $$
        H=a(x_m)b(x_1,\ldots,x_{m-1})y+f(z,t)+x_m
        $$ such that
        $b \neq 0$ and $a(x_m)$ has at least two distinct roots in $\overline{k}$. Then the following are equivalent:
        \begin{enumerate}[\rm (a)]
            \item $\mathrm{Spec}\,A$ is smooth and $\A^1$-contractible.

            \item $A=k^{[m+2]}$.

            \item $k[z,t]=k[f]^{[1]}$.

            \item $k[x_1,\ldots,x_m,y,z,t]=k[H]^{[m+2]}$.
        \end{enumerate}
    \end{thm}


     \begin{thm}
         Let $k$ be any field
        $$A=\frac{k[x_1,\ldots,x_m,y,z,t]}{\left( x_1^{r_1}a_1(x_1) \cdots x_{m-1}^{r_{m-1}}a_{m-1}(x_{m-1})x_m^{r_m}a_m(x_m) y + f(z, t) + x_m \right)}, \quad r_i>1,
       $$        
and 
\begin{equation}\label{b2}
H=x_1^{r_1}a_1(x_1) \cdots x_{m-1}^{r_{m-1}}a_{m-1}(x_{m-1})x_m^{r_m}a_m(x_m) y + f(z, t) + x_m, \ a_i(x_i) \neq 0.
\end{equation}
Then the following statements are equivalent:
		\begin{enumerate}[\rm (a)]
		\item $A=k^{[m+2]}$ 
		\item $k[Z,T]=k[f]^{[1]}$
		\item $k[X_1,\ldots,X_m,Y,Z,T]=k[H]^{[m+2]}$.
		\end{enumerate}
	 \end{thm}
   Note that the above theorems identify a large subfamily of varieties of the form $\mathrm{Spec}\,A$, as in \eqref{a}, which are actually isomorphic to the affine space ($\A_k^{m+2}$), and  in fact the equivalences of \thref{embedding1}(a), (b), (d) and \thref{embedding}(a) and (c) determines that the embeddings of the corresponding hyperplane into the ambient space $\A_k^{m+3}$ are indeed rectifiable, and thus producing new evidences towards the Abhyankar-Sathaye conjecture. 

   In the last section, in every dimension $\geqslant 4$ we construct examples of smooth exotic $\mathbb{A}^1$-contractible varieties which are isomorphic to none of the known examples in \cite{DG}.

    \subsection*{Organization of the article}
    The paper is organized as follows. In Section~2, we address the necessary part of Question~1.1. In Subsection~2.1 provides preliminary results that will be subsequently used in this paper, Subsection~2.2 develops the required technical machinery from $G$-theory and higher Chow groups, and Subsection~2.3 contains the main necessary criterion (Theorem 1.2) for $\A^1$-contractibility of the varieties $\mathrm{Spec}\,A$ defined by \eqref{a}; we also record there two applications of this criterion, namely that the factorial threefolds constructed in \cite{GGPJA} are not $\A^1$-contractible, and the construction of smooth affine varieties that are $\A^1$-chain connected but not $\A^1$-contractible. In Section~3, we turn to the sufficient part of Question~1.1. In Subsection~3.1 we study singular $\A^1$-contractible affine curves over perfect fields and prove Theorem~1.3. In Subsection~3.2 we establish the sufficient conditions for (stable) $\A^1$-contractiblity of the varieties $\mathrm{Spec}\,A$ (Theorems \ref{sufficientunstable}, \ref{sufficient}). In Section~4, we answer Question~1(ii) by establishing Theorems~1.4 and~1.5, and deduce the resulting rectifiability statements for the associated embeddings of affine spaces. Finally, in Section~5, we construct, in every dimension $\geqslant 4$, a new family $\mathcal{X}_a$ of smooth $\A^1$-contractible affine varieties (Theorem~\ref{0.2}); computing their Derksen and Makar-Limanov invariants (Propositions~\ref{0.3} and~\ref{0.4}), we show that they are not isomorphic to affine spaces (Corollary~\ref{0.5}) and are isomorphic to none of the previously known examples of \cite{DG} and \cite{Vij} (Proposition~\ref{0.6}).


	\section{Necessary conditions for $\A^1$-contractibility for Koras-Russell type varieties}

    Throughout this section, unless specified, 
    $H$ denotes a polynomial of the form
    \begin{equation}\label{2}
        H=a(x_m)b(x_1,\ldots,x_{m-1})y+f(z,t)+x_m, a,b \neq 0, \text{ and } f \notin k.
    \end{equation}
    We call the corresponding hypersurfaces {\it Koras--Russell type varieties}. 

    \subsection{Preliminary results.} We begin with some useful results from \cite{pre} that will be employed subsequently in this paper. The following result gives a criterion for a simple birational extension of a UFD to be a UFD.
    
    \begin{prop}\thlabel{ufdg} (\cite[Proposition 3.5]{pre})
		Let $R$ be a UFD. Let $u,v\in R\setminus\{0\}$ be such that $C=\frac{R[Y]}{(uY-v)}$ be an integral domain.
		We consider $R$ as a subring of $C$.
		Let $u:= \prod_{i=1}^{n}u_i^{r_i}$ be a prime factorization of $u$ in $R$. Suppose that for every $i \in \{ 1, \dots, n\}$ for which $(u_i,v)R$ is a proper ideal, we have $\prod_{j\neq i}^{}u_j^{s_j}\notin (u_i,v)R,$ for arbitrary integers $s_j\geqslant 0$. Then the following statements are equivalent:  
		\begin{enumerate}[\rm(i)]
			\item $C$ is a UFD.
			
			\item For each $i$, $1 \leqslant i \leqslant n$, either $u_i$ is prime in $C$ or $u_i \in C^*$.
			
			\item For each $i$, $1 \leqslant i \leqslant n$, either $(u_i,v)R$ is a prime ideal in $R$ or $(u_i,v)R = R$, i.e., the image of $v$ in $\frac{R}{u_iR}$ is either a prime in $\frac{R}{u_iR}$ or a unit in $\frac{R}{u_iR}$. 
		\end{enumerate}  
	\end{prop}

    \begin{rem}\thlabel{ufd}
        \rm 
       Let $A=\dfrac{k[x_1,\ldots,x_m,y,z,t]}{(a(x_m)b(x_1,\ldots,x_{m-1})y+f(z,t)+x_m)}$ be such that $A \otimes_k \overline{k}$ is a UFD and $(A \otimes_k \overline{k})^*=\overline{k}^*$. For every root $\lambda$ of $a(x_m)$ in $\overline{k}$, since $x_m-\lambda \notin (A \otimes_k \overline{k})^*$, by \thref{ufdg} ( (ii) $\Leftrightarrow$ (iii) ), it follows that $f(z,t)+\lambda$ is irreducible in $\overline{k}[z,t]$.
    \end{rem}

     The next result will play a crucial role in \thref{thm2}, which is one of the main steps to prove \thref{1.2}. 

        \begin{lem}\thlabel{lem4} (\cite[Lemma 4.6]{pre})
			Let $k$ be an algebraically closed field. 
			Suppose that $C$ is a regular affine $k$-domain and $R$ is a reduced affine $k$-algebra. If the map $R \hookrightarrow R\otimes_k C$ induces surjective maps $G_i(R) \rightarrow G_i(R\otimes_k C)$ for $i=0$ and $1$, then the canonical inclusion $\tau: k \hookrightarrow C$ induces isomorphisms of $K_i$-groups for $i=0$ and $1$ that means $K_0(C)=\bZ$ and $K_1(C)=k^*$ respectively.	
		\end{lem}

    The next result from \cite{pre} gives a criterion for a ring of the form
    \begin{equation}\label{b}
        A=\dfrac{k[x_1,\ldots,x_m,y,z,t]}{\left(  x_1^{r_1}a_1(x_1) \cdots x_{m-1}^{r_{m-1}}a_{m-1}(x_{m-1})x_m^{r_m}a_m(x_m) y + f(z, t) + x_m \right)}\,\text{~with~} r_i>1
    \end{equation}
    to be isomorphic to $k^{[m+2]}$. Before stating the result we recall the definition of exponential map and Makar-Limanov invariant of a $k$-algebra.

\smallskip
\noindent
{\bf Definition 2.}
Let $R$ be a $k$-algebra and let
$
\varphi_U : R \longrightarrow R[U]
$
be a $k$-algebra homomorphism. We say that $\varphi=\varphi_U$ is an
\emph{exponential map} on $R$ if $\varphi$ satisfies the following two properties:
\begin{enumerate}
    \item[(i)] $\varepsilon_0 \varphi_U$ is the identity on $R$, where
    $
    \varepsilon_0 : R[U] \longrightarrow R
    $
    is the evaluation map at $U=0$.

    \item[(ii)] $\varphi_V \varphi_U = \varphi_{V+U}$, where
    $
    \varphi_V : R \longrightarrow R[V]
    $
    is extended to a $k$-algebra homomorphism
    $
    \varphi_V : R[U] \longrightarrow R[U,V]
    $
    by defining $\varphi_V(U)=U$.
\end{enumerate}

Given an exponential map $\varphi$ on a $k$-algebra $R$, the ring of invariants of
$\varphi$ is a subring of $R$ given by
\[
R^\varphi
=
\{\, a \in R \mid \varphi(a)=a \,\}.
\]
We denote the set of all exponential maps on $R$ by
$\operatorname{EXP}(R)$.

The \emph{Makar-Limanov invariant} of $R$ is a subring of $R$ defined by
\[
\operatorname{ML}(R)
=
\bigcap_{\varphi \in \operatorname{EXP}(R)} R^\varphi.
\]

\begin{rem}
    If $R=k^{[n]}$, then $\ml(R)=k$.
\end{rem}

    \begin{thm}\thlabel{lin} (\cite[Theorem 5.22(ii)]{pre})
	Let $A$ be as in (\ref{b}) and $\ml(A)= k$. If $k$ is an infinite field then there exists a system of coordinates $\{z_1,t_1\}$ of $k[z,t]$ and $b_0,b_1 \in k^{[1]}$, such that $f(z,t)=b_0(z_1)+b_1(z_1)t_1$. 
	\end{thm}

	\begin{prop}\thlabel{thm1}
		Let $A$ be as in (\ref{a}). If $\mathrm{Spec}\, A$ is smooth $\mathbb{A}^1$-contractible, then for every root $\lambda$ of $a(x_m)$ in $\overline{k}$, $f(z,t)+\lambda$ is irreducible in $\overline{k}[z,t]$.
        Moreover, if $k$ is infinite, $A$ is of the form \eqref{b} and $A=k^{[m+2]}$, then up to a change of coordinates of $k[z,t]$, \( f(z, t) = b_0(z) + b_1(z) t \), such that \( \gcd(b_0(z), b_1(z)) = 1 \).
	\end{prop}
	\begin{proof}
		Since $\mathrm{Spec}\,A$ is smooth $\A^1$-contractible, $\mathrm{Spec}\,(A \otimes_k \overline{k})$ is so. Therefore, $A \otimes_k \overline{k}$ is a UFD and $(A \otimes_k \overline{k})^*=\overline{k}^*$. Therefore, $f(z,t)+\lambda$ is irreducible by \thref{ufd}.
		Moreover, if $A$ is as in \eqref{b} and $A=k^{[m+2]}$, then $\ml(A)=k$, and therefore the conclusion follows from \thref{lin}.
	\end{proof}

    \subsection{Technical Foundation: Higher Chow groups and G-theory.} We now establish the results that will be useful to prove \thref{1.2}. These results show how the $\mathbb{A}^1$-contractibility of the Koras--Russell type varieties $\Spec\, A$ (given by the rings \eqref{a}) translates into properties of Chow groups of the curves $\{f(z,t)+\lambda=0\}$, and they describe the possible singularities of these curves.

    \bigskip
	
	\begin{prop}\thlabel{prop1}
	     Let 
		$$
		A= \frac{k[x_1, \ldots, x_m, y, z, t]}
        { \left(  a(x_m)b(x_1,\ldots,x_{m-1}) y + f(z, t) + x_m \right) },\ b \neq 0 \text{ and } f \notin k
		$$ 
        If $\mathrm{Spec}\,A$ is smooth and $\mathbb{A}^1$-contractible,  then there exists a regular affine $\overline{k}$-domain $B$ such that $B \hookrightarrow B[z,t]/(f+\lambda)$ induces isomorphisms on $G_i$-groups and total higher Chow groups $\mathrm{CH}_*(-,i):=\bigoplus_{j \in \mathbb{Z}} \mathrm{CH}_j(-,i)$, $i \geqslant 0$ for every root $\lambda$ of $a(x_m)$ in $\overline{k}$.

	\end{prop}
\begin{proof}

         Since $\mathrm{Spec}\, A \to \mathrm{Spec}\,k$ is an $\mathbb{A}^1$-weak equivalence, for $\overline{A}:=A \otimes_k \overline{k}$, $\mathrm{Spec}\, \overline{A} \to \mathrm{Spec}\, \overline{k}$ is also an $\mathbb{A}^1$-weak equivalence. Therefore, without loss of generality, we assume that $k=\overline{k}$ and $A=\overline{A}$. 
         
		Note that 
		\[
		 R:=k[x_1, \ldots, x_m] \xhookrightarrow[]{\eta} A
		\]
		
		is a flat extension (cf. \cite[Corollary 2, 20.F]{Matc}). 
        Since both $\mathrm{Spec}\,R$ and $\mathrm{Spec}\,A$ are $\mathbb{A}^1$-contractible, $\eta$ induces an $\mathbb{A}^1$-weak equivalence and hence induces isomorphism of total higher chow groups and $G_i$-groups. Further $\overline{\eta}: R/p_1 \to A/p_1$ induces isomorphism of $G_i$ and total higher chow groups by $\mathbb{A}^1$-invariance. 
Let $p_1,\ldots p_n$ be prime divisors of $b(x_1,\ldots,x_{m-1})$. 

		Now, consider the commutative diagrams:
		\[
		\begin{tikzcd}[row sep=large, column sep=small]
			 \arrow{r} & G_{i+1}(R/p_1) \arrow{r} \arrow{d}{\cong} 
			& G_{i+1}(R) \arrow{r} \arrow{d}{\cong} 
			& G_{i+1}(R[p_1^{-1}]) \arrow{r} \arrow{d} 
			& G_i(R/p_1) \arrow{r} \arrow{d}{\cong} 
			& G_i(R) \arrow{d}{\cong} \\
			 \arrow{r}& G_{i+1}(A/p_1) \arrow{r}  & G_{i+1}(A)\arrow{r} & G_{i+1}(A[p_1^{-1}]) \arrow{r} &G_i(A/p_1) \arrow{r} &G_i(A).
		\end{tikzcd}
		\]
       and
       \[
		\begin{tikzcd}[row sep=large, column sep=small]
			 \arrow{r} & \mathrm{CH}_{*}(R/p_1,i+1) \arrow{r} \arrow{d}{\cong} 
			& \mathrm{CH}_{*}(R,i+1) \arrow{r} \arrow{d}{\cong} 
			& \mathrm{CH}_{*}(R[p_1^{-1}],i+1) \arrow{r} \arrow{d} 
			& \mathrm{CH}_{*}(R/p_1,i) \arrow{r} \arrow{d}{\cong} 
			& \mathrm{CH}_*(R,i) \arrow{d}{\cong} \\
			 \arrow{r}& \mathrm{CH}_{*}(A/p_1,i+1) \arrow{r}  & \mathrm{CH}_{*}(A,i+1)\arrow{r} & \mathrm{CH}_{*}(A[p_1^{-1}],i+1) \arrow{r} &\mathrm{CH}_*(A/p_1,i) \arrow{r} &\mathrm{CH}_*(A,i).
		\end{tikzcd}
		\]
        Therefore, by the Five Lemma, $R_1:=R[p_1^{-1}] \xhookrightarrow[]{\eta_1} A_1:=A[p_1^{-1}]$ induces isomorphisms on $G_i$-groups and higher Chow groups.

        Considering the commutative diagrams of localization sequences of $G_i$-groups and higher Chow groups induced by $\eta_1$ with respect to the element $p_2$, we obtain that $R_2:=R[p_1^{-1},p_2^{-1}] \xhookrightarrow[]{\eta_2} A_2:=A[p_1^{-1},p_2^{-1}]$ induces isomorphisms on $G_i$-groups and higher Chow groups respectively. Iterating this process for $p_3, \ldots,p_n$ successively, we obtain
         $$
        R_n:=R[p_1^{-1},\ldots,p_n^{-1}]=B[x_m] \xhookrightarrow[]{\eta_n} A_n:=A[p_1^{-1},\ldots,p_n^{-1}]=\dfrac{B[x_m,y,z,t]}{(a(x_m)b(x_1,\ldots,x_{m-1})y+f(z,t)+x_m)}
         $$
		induces isomorphisms on $G_i$-groups and higher Chow groups, where $B=k[x_1,\ldots,x_{m-1},p_1^{-1},\ldots,p_n^{-1}]$.

		Now let $q_1,\ldots,q_l$ be the prime divisors of $a(x_m)$ in $k[x_m]$. Then $\eta_n$ induces the following commutative diagram of localization sequences:

        \[
		\begin{tikzcd}[row sep=large, column sep=small]
			 \arrow{r} & G_{i}(R_n) \arrow{r} \arrow{d}{\cong} 
			& G_{i}(R_n[ \prod_{j=1}^lq_j^{-1}]) \arrow{r} \arrow{d}{\cong} 
			& G_{i-1}(R_n/(\prod_{j=1}^l q_j)) \arrow{r} \arrow{d} 
			& G_{i-1}(R_n) \arrow{r} \arrow{d}{\cong} 
			& G_{i-1}(R_n[ \prod_{j=1}^lq_j^{-1}]) \arrow{d}{\cong} \\
			 \arrow{r}& G_{i}(A_n) \arrow{r}  & G_{i}(A_n[\prod_{j=1}^l q_j^{-1}])\arrow{r} & G_{i-1}(A_n/(\prod_{j=1}^lq_j)) \arrow{r} &G_{i-1}(A_n) \arrow{r} &G_{i-1}(A_n[\prod_{j=1}^l q_j^{-1}]).
		\end{tikzcd}
		\]
       and
       
       \adjustbox{scale=0.9,center}{
       
		\begin{tikzcd}[row sep=large, column sep=tiny]
			 \arrow{r} & \mathrm{CH}_{*}(R_n,i) \arrow{r} \arrow{d}{\cong} 
			& \mathrm{CH}_{*}(R_n[q_1^{-1},\ldots,q_l^{-1}],i) \arrow{r} \arrow{d}{\cong} 
			& \mathrm{CH}_{*}\bigl(R_n/(q_1\cdots q_l),i-1\bigr) \arrow{r} \arrow{d} 
			& \mathrm{CH}_{*}(R_n,i-1) \arrow{r} \arrow{d}{\cong} 
			& \mathrm{CH}_*(R_n[q_1^{-1},\ldots,q_l^{-1}],i-1) \arrow{d}{\cong} \\
			 \arrow{r}& \mathrm{CH}_{*}(A_n,i) \arrow{r}  & \mathrm{CH}_{*}(A_n[q_1^{-1},\ldots,q_l^{-1}],i)\arrow{r} & \mathrm{CH}_{*}\bigl(A_n/(q_1\cdots q_l),i-1\bigr) \arrow{r} &\mathrm{CH}_*(A_n,i-1) \arrow{r} &\mathrm{CH}_*(A_n[q_1^{-1},\ldots,q_l^{-1}],i-1).
		\end{tikzcd}
		
       }
       
       \smallskip
        By the Five Lemma applied to the above diagrams, the middle arrows are isomorphisms. Since $q_j$'s are pairwise comaximal, we get that the inclusion $B \hookrightarrow B[z,t]/(f+\lambda)$ induces isomorphisms on $G_i$-groups and total higher Chow groups $\mathrm{CH}_*(-,i)$ for $i \geqslant 0$, for every root $\lambda$ of $a(x_m)$ in $\overline{k}$. 
        \end{proof}
    
The next result gives a criterion for an integral affine plane curve to have only unibranch singularities.

    \begin{prop}\thlabel{thm2}
        Let $k$ be an algebraically closed field, let $C=k[Z,T]/(g)$ be an integral domain, and let $B$ be a regular affine $k$-domain. If $B \hookrightarrow B\otimes C$ induces isomorphisms on $G_i$-groups and total higher Chow groups, then the normalization of $C$ is $k^{[1]}$ and $\mathrm{Spec}\, C$ has only unibranch singularities.     
    \end{prop}

\begin{proof}
     We first prove that the normalization of $C$ is $k^{[1]}$. Let $\widetilde{C}$ denote the normalization of $C$. Since $k=\overline{k}$ and $B$ is regular affine $k$-domain, it follows that the normalization $\widetilde{B\otimes_k C}$ of $B \otimes_k C$ is $B \otimes_k \widetilde{C}$. Let $\mathcal{I}$ be the conductor ideal of $\widetilde{C}$ into $C$ and $0 \neq a \in \mathcal{I}$. We now consider the following surjection
        \begin{equation}\label{a1}
        G_0(B \otimes_k C) \twoheadrightarrow G_0((B \otimes_k C) [1/a]).
        \end{equation}

        Now using the isomorphism $C[1/a] \xrightarrow{\cong} \widetilde{C}[1/a]$, it follows that $(B \otimes_k C)[1/a] \cong B \otimes_k \widetilde{C}[1/a]$.
        Since the inclusion $B \hookrightarrow B \otimes_k C$ induces isomorphism $G_0(B) \to G_0(B \otimes_k C)$, using the surjection \eqref{a1}, we get that $B \hookrightarrow B \otimes_k \widetilde{C}[1/a]$ induces a surjection 
        $$G_0(B) \twoheadrightarrow G_0(B \otimes_k \widetilde{C}[1/a]).$$ 
        We now note that $B$ is an affine $k$-domain.
        Therefore, by \thref{lem4}, we obtain $G_0(\widetilde{C}[1/a])= K_0(\widetilde{C}[1/a])=\mathbb{Z}$ so $\widetilde{C}[1/a]$ is a smooth rational curve and hence $C$ is rational and $\widetilde{C}\cong k[V, 1/h(V)]$, for some $h \in k^{[1]}$.  

        We now show that $\widetilde{C}^{*}=k^*$. For that we first consider the following natural maps:
        $$
        G_1(B) \xrightarrow{\phi_1} G_1(B \otimes_k\widetilde{C})=G_1(\widetilde{B \otimes_k C}) \xrightarrow{\phi_2} G_1(B \otimes_k C).
        $$
        Since $\phi_2 \phi_1$ is an isomorphism by our assumption, it follows that $\phi_1$ is injective. We now show that $\phi_1$ is surjective. For this, consider the following commutative square:

		\[
		\begin{tikzcd}[row sep=large, column sep=small]
			  G_1(\widetilde{B \otimes_k C})=G_{1}(B \otimes_k \widetilde{C}) \arrow[r,"\phi_2"] \arrow{d} 
			& G_{1}(B \otimes_k C) \arrow{d}\\
			  G_{1}((B \otimes_k \widetilde{C})[1/a]) \arrow{r}{\cong}  & G_{1}((B \otimes_k C)[1/a]).
		\end{tikzcd}
		\]
        Since $k$ is algebraically closed, and $\widetilde{C}=k[V,1/h(V)]$ it follows that the left vertical arrow is injective. Therefore, $\phi_2$ must be injective. Since $\phi_2\phi_1$ is isomorphism, it follows that $\phi_1$ is surjective and hence $\phi_1$ is an isomorphism. Therefore, by \thref{lem4}, it follows that $\widetilde{C}^*=k^*$. Therefore $\widetilde{C}=k^{[1]}$.

 We now show that $C$ has only unibranch singularities. Let us consider the normalization map $\eta: C \hookrightarrow k[u]:=\widetilde{C}$ and let $I=(p(u))\widetilde{C}$ be the conductor ideal. As $B \hookrightarrow B \otimes_k C$ induces isomorphisms on higher Chow groups, the normalization $\eta_B: B \otimes_k C \hookrightarrow \widetilde{(B \otimes_k C)}=B \otimes_k \widetilde{C}=B[u]$
 induces isomorphisms on the higher Chow groups as well. Moreover, $\eta_B$ induces an isomorphism of rings $(B \otimes_k C)[p(u)^{-1}] \xrightarrow{\cong} \widetilde{(B \otimes_k C)}[p(u)^{-1}]$. 
 Let $R:=B \otimes_k C$ and $\widetilde{R}:=\widetilde{(B \otimes_k C)}=B \otimes_k \widetilde{C}$.
 Since $\eta_B$ is a finite map, it induces the following commutative diagram of higher Chow groups

		\begin{tikzcd}[row sep=large, column sep=tiny]
			 \arrow{r} & \mathrm{CH}_{j}(\widetilde{R},i+1) \arrow{r} \arrow{d}{\cong} 
			& \mathrm{CH}_{j}(\widetilde{R}[p(u)^{-1}],i+1) \arrow{r} \arrow{d}{\cong} 
			& \mathrm{CH}_{j}(\widetilde{R}/(p(u)),i) \arrow{r} \arrow{d} 
			& \mathrm{CH}_{j}(\widetilde{R},i) \arrow{r} \arrow{d}{\cong} 
			& \mathrm{CH}_j(\widetilde{R}[ (p(u))^{-1}],i) \arrow{d}{\cong} \\
			 \arrow{r}& \mathrm{CH}_{j}(R,i+1) \arrow{r}  & \mathrm{CH}_{j}(R[p(u)^{-1}],i+1)\arrow{r} & \mathrm{CH}_{j}(R/(p(u)),i) \arrow{r} &\mathrm{CH}_j(R,i) \arrow{r} &\mathrm{CH}_j(R[(p(u))^{-1},i]).
		\end{tikzcd}
		
    Therefore, $\eta_B$ induces isomorphism 
    $$\mathrm{CH}_j(\widetilde{R}/p(u)) \xrightarrow{\cong} \mathrm{CH}_j(R/p(u)).$$ 
    Now note that the support of $\mathrm{Spec}(C/(p(u)))$ gives the singular locus of $\mathrm{Spec}\, C$. 
    Now if we consider the map  $\eta^{*}: \mathbb{A}_k^1 \to \mathrm{Spec}\, (C)$,
    from the above isomorphism of Chow groups for $j=\mathrm{dim}(B)$, it follows that the fiber over every singular point of $\mathrm{Spec}\, C$ contains exactly one point.  
    Therefore it follows that
    $\eta^{*}: \mathbb{A}_k^1 \to \mathrm{Spec}\, (C)$ is a bijection. Hence all the singularities of $\mathrm{Spec}\, C$ must be unibranched. 

\end{proof}



    \subsection{Necessary condition for $\A^1$-contractibility.} 
    
    We now prove \thref{1.2}, which gives the necessary conditions for $\mathbb{A}^1$-contractibility for Koras-Russell type varieties given by \eqref{a}.

\begin{thm}\thlabel{thm3}

             Let 	$A$ be a ring as in as in \eqref{a},   $\lambda$ be a root of $a(x_m)$ in $\overline{k}$, $\Gamma_{\lambda}:=\mathrm{Spec}\,(k(\lambda)[z,t]/(f+\lambda))$ and $\overline{\Gamma}_{\lambda}:=\Gamma_{\lambda} \times_{k(\lambda)} \overline{k}$.  
            If $\mathrm{Spec}\, A$ is a smooth $\mathbb{A}^1$-contractible affine variety. Then

            \begin{itemize}
                \item[(a)]
             $\Gamma_{\lambda}$ has only unibranch singularities and the normalization of $\overline{\Gamma}_{\lambda}$ is isomorphic to $\mathbb{A}_{\overline{k}}^1$.  Furthermore, the corresponding Nisnevich sheaves  $h_{\overline{\Gamma}_{\lambda}}, h_{\mathbb{A}_{\overline{k}}^1}$ over $Sm_{\overline{k}}$ are isomorphic. 
             
                \item[(b)] If $k$ is a perfect field, then the curve $\Gamma_{\lambda}$ is a polynomial curve; that is, the normalization of $\Gamma_{\lambda}$ is isomorphic to $\mathbb{A}_{k(\lambda)}^1$. Furthermore, the corresponding Nisnevich sheaves $h_{\Gamma_{\lambda}}, h_{\mathbb{A}_{k(\lambda)}^1}$ over $Sm_{k(\lambda)}$ are isomorphic.
            \end{itemize}

	\end{thm}

	\begin{proof}
     (a) Since $\mathrm{Spec}\, A$ is $\A^1$-contractible, $\mathrm{Spec}\, (A \otimes_k \overline{k})$ is $\A^1$-contractible, therefore $A$ is a geometrically factorial domain and hence $\overline{k}[Z,T]/(f+\lambda)$ is an integral domain (see \thref{thm1}). 

     \smallskip
     \noindent
     Now by \thref{prop1}, there exists a regular $\overline{k}$-domain $B$ such that $B \hookrightarrow B[Z,T]/(f+\lambda)$ induces isomorphisms on higher Chow groups and $G_i$-groups. Now by \thref{thm2}, the normalization of $\overline{\Gamma}_{\lambda}$ is $\mathbb{A}_{\overline{k}}^1$ and $\overline{\Gamma}_{\lambda}:=\Gamma_{\lambda} \times_{k(\lambda)} \overline{k}$ has only unibranch singularities. Therefore $\Gamma_{\lambda}$ also has only unibranch singularities. 
     Since $\overline{\Gamma}_{\lambda}$ has only unibranch singularities, the normalization map $\mathbb{A}_{\overline{k}}^1 \to \overline{\Gamma}_{\lambda}$ is bijective. Therefore, the corresponding Nisnevich sheaves $h_{\overline{\Gamma}_{\lambda}}$ and $h_{\mathbb{A}_{\overline{k}}^1}$ on $Sm_{\overline{k}}$ are isomorphic  (cf. \cite[Example 2.1]{AD}).

     \smallskip
     \noindent
     (b)
     Let $\widetilde{\Gamma}_{\lambda}$ be the normalization of $\Gamma_{\lambda}$. If $k$ is a perfect field, the base change $\widetilde{\Gamma}_{\lambda} \times_{k(\lambda)} \overline{k}$ is the normalization of $\overline{\Gamma}_{\lambda}$. 
     Therefore, $\widetilde{\Gamma}_{\lambda}\times_{k(\lambda)} \overline{k} \cong \mathbb{A}_{\overline{k}}^1$. Since $k$ is perfect, it follows that $\widetilde{\Gamma}_{\lambda} \cong \mathbb{A}_{k(\lambda)}^1$. Since the map $\A_{\overline{k}}^1 \to \overline{\Gamma}_{\lambda}$ is bijective, the normalization map $\mathbb{A}_{k(\lambda)}^1 \to \Gamma_{\lambda}$ is universally bijective and hence the map of Nisnevich sheaves $h_{\A_{k(\lambda)}^1} \to h_{\Gamma_{\lambda}}$ on $Sm_{k(\lambda)}$ is an isomorphism (cf. \cite[Example 2.1]{AD}).
\end{proof}

In \cite{GGPJA}, examples of smooth affine threefolds that are factorial but not isomorphic to affine space are constructed over arbitrary fields. Using \thref{thm3}, we obtain the stronger conclusion below that these varieties cannot even be $\A^1$-contractible.

\noindent
\textbf{Example 1.}
Let \(r,s>1\) be relatively prime integers, and let \(a(x)\in k[x]\)
have at least two distinct roots in \(\overline{k}\) and let
\[
\mathcal{X}=\operatorname{Spec}\left(
\frac{k[x,y,z,t]}{(a(x)y+z^r+t^s+x)}\right),
\]
be a smooth affine variety.
Note that $\sigma: \mathcal{X} \to \mathrm{Spec}\,k[x,z,t]$ is a simple birational extension
with respect to the divisor \(\mathcal{D}:=V(a(x))\), center \(\mathcal{Z}:=V(a(x),z^r+t^s+x)\). 
Note that
 $\mathcal{E}:=V(a(x),z^r+t^s+x) \times \mathbb{A}^1$ is the exceptional divisor of the extension. 
 Over \(\mathbb C\),  \(\mathcal{X}\)
is not topologically contractible. 
If \(\mathcal X(\mathbb C)\) were contractible, then \(\sigma\) would induce isomorphisms on singular homology, and \cite[Theorem~ 3.1]{afmod} would imply that \(\sigma|_{\mathcal E}\) induces homology isomorphisms.
 But this is a contradiction as $a(x)$
has at least two distinct roots and hence it has a non-zero root $\lambda$ and the restriction of $\sigma$ to the component $\mathcal{E}_{\lambda}:= V(z^r+t^s+\lambda) \times \mathbb{A}^1$,
given by $\sigma|_{\mathcal{E}_{\lambda}}: \mathcal{E}_{\lambda} \to \mathrm{Spec}\, \mathbb{C}[z,t]$
does not induce an isomorphism on singular homology.

Over any field of characteristic zero, \(\mathcal{X} \ncong\mathbb A^3_k\) (\cite{KVZ}, \cite{KO}, \cite{Mau})
and  if every root of \(a(x)\) in \(\overline{k}\) is multiple, then over a field of characteristic $ p>0$, 
\(\mathcal{X}\ncong\mathbb A^3_k\) by \cite[Theorem~4.6]{GGPJA}. 
However, Theorem~2.9(a) gives a stronger conclusion that $\mathcal{X}$ is not even \(\mathbb A^1\)-contractible over any field of characteristic zero
and over a field of positive characteristic $p>0$ such that $p\nmid rs$,
because for a non-zero root $\lambda$ of $a(x)$, the curve $z^r+t^s+\lambda=0$ in $\mathbb A_{\overline{k}}^2$ is smooth and not isomorphic to $\mathbb A_{\overline{k}}^1$; in particular, it is not a polynomial curve. \qed

We now construct an example of a smooth $\A^1$-chain connected affine variety that is not $\A^1$-contractible.
 Let \(k\) be any field, let \(m\ge2\), and let
\(f\in k[z,t] \setminus k\). Let
\begin{equation}\label{a.1}
C:=V(f)\subset\mathbb A^2_k,
\qquad
\mathcal{X}_{m,f}=\operatorname{Spec}\left(
\frac{k[x_1,\ldots,x_m,y,z,t]}{(x_1^{r_1}\ldots x_m^{r_m} y+f(z,t)+x_m)}\right).
\end{equation}
The following result gives a criterion for $\mathbb{A}^1$-chain connectedness of $\mathcal{X}_{m,f}$.

\begin{prop}\label{0.1}
     Let \(\mathcal{X}_{m,f}\) be a smooth affine variety as in \eqref{a.1} and let
$C$ be as in \eqref{a.1} such that it has a smooth \(k\)-rational point. 
Then \(\mathcal{X}_{m,f}\) is \(\mathbb A^1\)-chain connected.   
\end{prop}

\begin{proof}
We prove this by induction on $m$. Let $m=1$ and 
$$
\mathcal{X}_{1,f}=\operatorname{Spec}\left(
\frac{k[x_1,y,z,t]}{(x_1^{r_1}y+f(z,t)+x_1)}\right).
$$
Let $C$ be smooth at \(c=(z_0,t_0)\in C(k)\). 
Let $L$ be a finitely generated field extension of $k$.
Since $C_L:=C \times_k L$ is smooth at \(c\), the Jacobian criterion gives \(f_z(z_0,t_0)\ne0\) or \(f_t(z_0,t_0)\ne0\).
Therefore, after interchanging \(z\) and \(t\) if
necessary, we may assume \(f_t(z_0,t_0)\ne0\).
We now construct a polynomial $T(u) \in L[u]$ such that \(T(0)=t_0\) and \(f(z_0,T(u))+u\equiv0\pmod{u^{r_1+1}}\).
Let \(T_0(u)=t_0\), and define \(T_j(u)\) recursively for \(j\geqslant 0\) by
 \[
T_{j+1}(u)=T_j(u)-\frac{e_{j+1}}{f_t(z_0,t_0)}u^{j+1}
\]
where \(e_{j+1}\) is the coefficient of \(u^{j+1}\) in \(f(z_0,T_j(u))+u\).
We define $T(u):=T_{r_1}(u)$. By induction we show that 
$$
f(z_0,T_{\widetilde{r}}(u))+u \equiv 0 \pmod{u^{\widetilde{r}+1}} \quad \text{for every~} \tilde{r} \in \mathbb{N}.
$$
We first note that
$f(z_0,T_1(u))+u \equiv 0 \pmod{u^2}$. Suppose that \(f(z_0,T_j(u))+u \equiv 0 \pmod{u^{j+1}}\) for \(j<\tilde{r}\). 
We now note that $$
f(z_0, T_{j+1}(u))+u \equiv f(z_0, T_j(u))+u -e_{j+1} u^{j+1} \pmod{u^{j+2}}.
$$
Therefore, by induction hypothesis, we have \(f(z_0,T_{j+1}(u))+u \equiv 0 \pmod{u^{j+2}}\),hence the induction step is done.

We now define a polynomial \(Y(u)\) by
\[
Y(u):=-\frac{f(z_0,T(u))+u}{u^{r_1}}\in L[u].
\]
We now consider the projection map $p_{x_1}: \mathcal{X}_{1,f} \to \mathbb{A}_k^1$. 
Let $\mu \in L^*$. Then
$
\gamma_{\mu}:\mathbb A^1_L\longrightarrow X_{1,f}$ defined by
$$
\gamma_{\mu}(u):=\bigl(\mu u,Y(\mu u),z_0,T(\mu u)\bigr)
$$ joins the point \(\gamma_{\mu}(0)=(0,Y(0),z_0,t_0)=(0,0,z_0,t_0)\) from $p_{x_1}^{-1}(0)$ to the point \(\gamma_{\mu}(1)=(\mu,Y(\mu),z_0,T(\mu))\) on $p_{x_1}^{-1}(\mu)$.
Since every fiber $p_{x_1}^{-1}(\mu)$ for $\mu \in L^*$ is isomorphic to $\mathbb{A}^2_L$, it is $\mathbb{A}^1$-chain connected.
Therefore, any point on $p_{x_1}^{-1}(\mu)$ is joined to $(0,0,z_0,t_0)$ by a chain of copies of $\mathbb{A}^1_L$. 
Now a point $(0,y_0,z_0,t_0)$ on $p_x^{-1}(0)(L)$ 
is joined to $(0,0,z_0,t_0)$ by a map $\gamma_{1}: \mathbb{A}_L^1 \to \mathcal{X}_{1,f}$ given by
$$\gamma_{1}(u):=\bigl(0, (1-u)y_0,z_0,t_0\bigr).$$
Now any $L$-rational point on the special fiber $p_x^{-1}(0)$ is of the form \((0,y_0,z_*,t_*)\) for some \(y_0\in L\) and \((z_*,t_*)\in C(L)\).
Note that there is a closed embedding
\[
\iota:\mathbb A^2_k\hookrightarrow X_{1,f} \text{~given by~}
(z,t)\longmapsto\bigl(-f(z,t),0,z,t\bigr).
\]
So the point \((0,0,z_*,t_*)\) is joined to \((0,0,z_0,t_0)\) by a chain of copies of \(\mathbb A^1_L\).
Further a point \((0,y_0,z_*,t_*)\) is joined to \((0,0,z_*,t_*)\) by image of a line $\mathbb{A}^1_L$ given by $\gamma_2:\mathbb{A}^1_{L} \to \mathcal{X}_{1,f}$ such that $$\gamma_2(u):=(0, (1-u)y_0,z_*,t_*).$$
Thus any point on the special fiber \(p_x^{-1}(0)(L)\) is joined to \((0,0,z_0,t_0)\) by a chain of copies of \(\mathbb A^1_L\), and hence $\mathcal{X}_{1,f}$ is $\A^1$-chain connected. 

Let us now assume $m>1$ and $\mathcal{X}_{m-1,f}$ is  $\mathbb{A}^1$-chain connected. 
We now show that $\mathcal{X}_{m,f}$ is $\mathbb{A}^1$-chain connected.

Now consider the projection
\[
\widetilde p: \mathcal{X}_{m,f}\longrightarrow\Spec k[x_1].
\]
Note that
$
\widetilde p^{-1}(0)\cong\A_k^{m+1},
$
which is $\A^1$-chain connected.

For $\mu\in L^*$, we have
\[
\widetilde p^{-1}(\mu)
\cong
\Spec\left(
\frac{k(\mu)[x_2,\ldots,x_m,y,z,t]}
{\bigl(x_2^{r_2}\cdots x_m^{r_m} y+f(z,t)+x_m\bigr)}
\right) \cong (\mathcal{X}_{m-1,f})_{k(\mu)}.
\]
This is $\A^1$-chain connected by the induction hypothesis.
Now consider the path
\[
\begin{aligned}
\gamma_\mu:\A_L^1&\longrightarrow \mathcal{X}_{m,f},\\
u&\longmapsto
\bigl(\mu u,\underbrace{1,\ldots,1}_{m-2},\mu u,
\widetilde{Y}(\mu u),z_0,T_{r_1+r_m}(\mu u)\bigr),
\end{aligned}
\]
where
\[
\widetilde{Y}(\mu u):=-\frac{f(z_0,T_{r_1+r_m}(\mu u))+\mu u}{(\mu u)^{r_1+r_m}}\in L[u],
\]
 Let $P=\gamma_{\mu}(1)$, for some $\mu \in L^*$. Then any point in the fiber $\widetilde p^{-1}(\mu)$ is joined to $P$
by a chain of $\mathbb{A}^1_L$, and then
to $q_0:=\gamma_\mu(0)$ by the path
$\gamma_\mu$. Any point in $\widetilde p^{-1}(0)$
is joined to $q_0$ by a chain of $\A_L^1$. Thus any two $L$-point of
$\mathcal{X}_{m,f}$ are joined by a chain of $\mathbb{A}^1_L$.
Thus we have proved that $\mathcal{X}_{m,f}$ is $\A^1$-chain connected.
\end{proof}

\smallskip
\noindent
{\bf Example 2.}
Let 
$$
\mathcal{X}_{m,f}:= \mathrm{Spec}\left(\dfrac{k[x_1,\ldots,x_m,y,z,t]}{(x_1^{r_1}\cdots x_m^{r_m}y+z^2+t^3-t^2+x_m)}\right), \quad r_i>1
$$ 
where $f(z,t)=z^2+t^3-t^2$ defines a curve in $\mathbb{A}^2_k$ with a node at the origin.
By the above Proposition \ref{0.1}, $\mathcal{X}_{m,f}$ is $\mathbb{A}^1$-chain connected.
However, it is not \(\mathbb A^1\)-contractible over any field of characteristic \( \neq 2\) by Theorem 2.9(a). \qed

    \section{Sufficient condition for $\A^1$-contractibility}
    
    In this section we prove sufficient conditions for Koras--Russell type varieties to be stably $\A^1$-contractible. We begin with a general result on the classification of singularities of $\mathbb{A}^1$-contractible affine curves.

    \subsection{$\A^1$-contractible Curves and their singularities.} 


      

 The next result classifies the singularities of a geometrically integral $\A^1$-contractible affine curve over a perfect field.

\begin{thm}\thlabel{uni-branched}
Let \(k\) be a perfect field, and let
\(X\in Sch_k\) be a
geometrically integral affine \(k\)-curve. If \(X\) is
\(\A^1\)-contractible in $\mathcal{H}((Sm_k)_{Nis},\A^1)$, then $X$ has only unibranch singularities.
\end{thm}

\begin{proof}
Let $\overline{X}= X \times_k \overline{k}$. Since $X$ is $\A^1$-contractible in $\mathcal{H}((Sm_k)_{Nis},\A^1)$,  $\overline{X}$ is also $\A^1$-contractible in $\mathcal{H}((Sm_{\overline{k}})_{Nis}, \A^1)$. 
Further, if $\overline{X}$ has only unibranch singularities, then so does $X$. 
Therefore, without loss of generality we assume $k$ is algebraically closed and $X=\overline{X}$.  
Let
\[
\nu\colon \widetilde X\longrightarrow X
\]
be the normalization of $X$. Let $U\subset X$ be the largest open subset over
which $\nu$ is an isomorphism and let
$
Z=(X\setminus U)_{red}.
$
Since $X$ is a noetherian integral curve, $Z$ is a finite reduced
$k$-scheme. Let
\[
Z'=\widetilde X\times_X Z,
\qquad
E=(Z')_{red}.
\]
We obtain a cartesian square
\begin{equation}
\label{eq:abstract-blowup}
\begin{tikzcd}
Z' \ar[r] \ar[d] & \widetilde X \ar[d,"\nu"]\\
Z \ar[r,"i"] & X .
\end{tikzcd}
\end{equation}
The morphism $\nu$ is finite, hence proper, and it induces an isomorphism
\[
\widetilde X\setminus Z'
\xrightarrow{\;\sim\;}
X\setminus Z.
\]
Thus \eqref{eq:abstract-blowup} is an abstract blow-up square.
Homotopy $K$-theory satisfies descent for abstract blow-up squares; (\cite[Theorem~C]{KST}, \cite[Theorem~3.9]{Cisinski2013}). Hence it induces
\[
\begin{tikzcd}
{KH}(X) \ar[r] \ar[d]
  & KH(\widetilde X) \ar[d]\\
KH(Z) \ar[r]
  & KH(Z')
\end{tikzcd}
\]
which is homotopy cartesian. Hence
\begin{equation}
\label{eq:KH-cdh}
{KH}(X)
\simeq
{KH}(\widetilde X)
\times_{\mathrm{KH}(Z')}
{KH}(Z).
\end{equation}
Since the homotopy $K$-theory is nilpotent invariant (\cite[Corollary~IV.12.5]{weibel}),
\[
{KH}(Z')\simeq {KH}(E).
\]
and hence
\begin{equation}
\label{eq:KH-cdh-reduced}
{KH}(X)
\simeq
{KH}(\widetilde X)
\times_{{KH}(E)}
{KH}(Z).
\end{equation}
Note that
$
\widetilde X, E, Z
$
are smooth $k$-schemes. 

We now show that 
\begin{equation}
\label{eq:sheaf-pushout}
\begin{tikzcd}
h_E \ar[r] \ar[d] & h_{\widetilde X} \ar[d]\\
h_Z \ar[r] & h_X
\end{tikzcd}
\end{equation}
 is a pushout square in $(Sm_k)_{Nis}$.   
Let $P$ be the pushout of the diagram $h_Z \longleftarrow h_E \longrightarrow h_{\widetilde X}$ in the category of Nisnevich presheaves on $Sm_k$. So
There is a canonical morphism of presheaves
\begin{equation}\label{eq:phi}
\phi:P\longrightarrow h_X.
\end{equation}
Now the pushout of the diagram
$h_{Z}\longleftarrow h_E\longrightarrow h_{\widetilde X}$ in the category of Nisnevich sheaves on $Sm_k$ is given by the Nisnevich sheafification $P^{Nis}$ of $P$. Since $h_X$ is a sheaf on $(Sm_k)_{Nis}$, $\phi$ induces a morphism of sheaves 
\begin{equation}\label{eq:sheafiso}
\widetilde{\phi}: P^{Nis} \longrightarrow h_X.
\end{equation}
We now show that $\widetilde{\phi}$ is an isomorphism of sheaves.
We first show that for every connected $V\in Sm_k$, the map
\begin{equation}\label{eq:connected-bij}
\phi(V):P(V)\longrightarrow h_X(V)
\end{equation}
is bijective. 

Let $g:V\to X$. Since $V$ is smooth and connected and $X$ is an integral curve, either $g$ is dominant or its image a closed point.
If $g$ is dominant, then the universal property of normalization gives a unique factorization
\[
V\xlongrightarrow{\widetilde g}\widetilde X\xrightarrow{\nu}X.
\]
If $g(V)=\{p\}$ is a point, then there are two cases. If $p\in Z$, then, $g$ factors through $Z$. If $p\notin Z$, then $p\in X \setminus Z$ which is isomorphic to $\widetilde X \setminus Z^{\prime}$. Therefore, $g$ factors through $\widetilde X$.
 Therefore every element of $h_X(V)$ is in the image of $\phi(V)$.

We now show the injectivity of $\phi(V)$. 
Suppose two 
\[
f,g \in h_{\widetilde X}(V)\amalg_{h_E(V)} h_Z(V)
\]
have the same image in $h_X(V)$. 
If $f,g \in h_Z(V)$, then they are equal because $Z\hookrightarrow X$ is a monomorphism.

Suppose $f:V\to\widetilde X$ and $g:V\to Z$, with same images in $h_X(V)$ that means,
\[
\nu\circ f=i\circ g.
\]
By the universal property of the fibre product $Z'=\widetilde X\times_X Z$, there is a unique morphism $\alpha: V\to Z'$ which induces $f$ and $g$. Since $V$ is reduced, every morphism $V\to Z'$ factors uniquely through $E=(Z')_{\mathrm{red}}$. Hence $f$ and $g$ are identified in the pushout $P(V)$.

Now suppose $f,g:V\to\widetilde X$ satisfy
\[
\nu f=\nu g.
\]
Let $\alpha:= \nu f=\nu g$.
If $\alpha$ is dominant, the universal property of normalization gives $f=g$. If the image of $\alpha$ is a point of $X\setminus Z$, then $f=g$ as $\nu: {\widetilde X \setminus Z^{\prime}} \to X \setminus Z$ is an isomorphism. If the image of $\alpha$ lies in $Z$, then $\alpha$ factors through some $\beta:V\to Z$, and each pair $(f,\beta)$ and $(g, \beta)$ determines a morphism $V \to Z^{\prime}$, and as $V$ is reduced a morphism $V\to E$. Thus both $f$ and $g$ represent  the same element in $P(V)$. Thus \eqref{eq:connected-bij} is injective.

Now let $V\in Sm_k$ be arbitrary. Since $V$ is of finite type over $k$, it has finitely many connected components
\[
V=\coprod_{j=1}^r V_j.
\]
 Hence $\{V_j\to V\}$ is a finite Zariski covering, and therefore a Nisnevich covering. By the preceding paragraph, $\phi(V_j)$ is bijective for every $j$.
Therefore by the Stacks Project, \cite[Lemma 7.10.17, Tag 0H6Y ]{stackproject}, the induced map on sheafifications
$$
\widetilde {\phi}: P^{Nis}\longrightarrow h_X
$$
is an isomorphism. Therefore \eqref{eq:sheaf-pushout} is a pushout square in $(Sm_k)_{Nis}$.
The morphism
\[
h_E\longrightarrow h_{\widetilde X}
\]
is a monomorphism. By \cite[Theorem~2.3.2]{MV}, monomorphisms are cofibrations in the
$\A^1$-local injective model structure on simplicial Nisnevich sheaves; in particular, all
objects are cofibrant. Hence \cite[Proposition~A.2.4.4]{LurieHTT} implies that the strict
pushout square \eqref{eq:sheaf-pushout} is homotopy pushout in
$\mathcal{H}((Sm_k)_{Nis},\mathbb{A}^1)$. Therefore, after $\mathbb P^1$-stabilization we obtain a homotopy pushout in
$SH(k)$:
\begin{equation}
\label{eq:stable-pushout}
\begin{tikzcd}
\Sigma^\infty_{\mathbb{P}^1} E_+ \ar[r] \ar[d]
&
\Sigma^\infty_{\mathbb{P}^1}\widetilde X_+ \ar[d]
\\
\Sigma^\infty_{\mathbb{P}^1} Z_+ \ar[r]
&
\Sigma^\infty_{\mathbb{P}^1}(h_X)_+ .
\end{tikzcd}
\end{equation}
Now applying the derived mapping-spectrum
functor
\[
\mathrm{RMap}_{SH(k)}(-,KGL_k)
\]
to \eqref{eq:stable-pushout} we get
\begin{align}\label{eq:mapping-pullback}
&\mathrm{RMap}_{SH(k)}
\left(
\Sigma^\infty_{\mathbb{P}^1}(h_X)_+,KGL_k
\right)\\
\nonumber
&\simeq
\mathrm{RMap}_{SH(k)}
\left(
\Sigma^\infty_{\mathbb{P}^1}\widetilde X_+, KGL_k
\right)
\times^h_{
\mathrm{RMap}_{SH(k)}
\left(
\Sigma^\infty_{\mathbb{P}^1}E_+,KGL_k
\right)}
\mathrm{RMap}_{SH(k)}
\left(
\Sigma^\infty_{\mathbb{P}^1}Z_+,KGL_k
\right).
\end{align}
Since $\widetilde X$, $E$, and $Z$ are smooth over $k$, the motivic
spectrum $KGL_k$ represents homotopy invariant $K$-theory on each of them (\cite[Theorem~2.20]{Cisinski2013}). Thus
from \eqref{eq:KH-cdh-reduced} and \eqref{eq:mapping-pullback} we get
\begin{equation}
\label{eq:singular-representability}
{KH}(X)
\simeq
\mathrm{RMap}_{SH(k)}
\left(
\Sigma^\infty_{\mathbb{P}^1}(h_X)_+,KGL_k
\right).
\end{equation}
Since
\[
h_X\longrightarrow h_{\Spec k}
\]
is a Nisnevich $\A^1$-weak equivalence, we have
\[
\Sigma^\infty_{\mathbb{P}^1}(h_X)_+
\simeq
\Sigma^\infty_{\mathbb{P}^1}(\Spec k)_+
=
\mathbf 1_k
\qquad\text{in } SH(k).
\]
Using \eqref{eq:singular-representability}, we get that
\[
\begin{aligned}
{KH}(X)
&\simeq
\mathrm{RMap}_{SH(k)}
\left(
\Sigma^\infty_{\mathbb{P}^1}(h_X)_+, KGL_k
\right)\\
&\simeq
\mathrm{RMap}_{SH(k)}
\left(
\mathbf 1_k, KGL_k
\right)\\
&\simeq
{KH}(\Spec k).
\end{aligned}
\]
Therefore
\[
KH_n(X)\cong KH_n(\Spec k)
\qquad\text{for every }n\in\mathbb Z.
\]
Since $\dim X=1$, 
\[
K_{-1}(X)\cong KH_{-1}(X),
\]
 this follows from
\cite[Theorem~B(ii)]{KST}. Therefore
\[
K_{-1}(X)\cong KH_{-1}(\Spec k)=0.
\]
Let
\(X=\operatorname{Spec}R\), let \(\widetilde R\) be the normalization of
\(R\), and let
$
  I=\operatorname{Ann}_R(\widetilde R/R)
$
be the conductor ideal which is an ideal of both \(R\) and \(\widetilde R\).
  For a ring \(A\), let
$
  h_0(A)=\#\pi_0(\operatorname{Spec}A).
$
By the Bass--Murthy theorem
\[
  K_{-1}(R)\cong \mathbb Z^{\rho},
  \qquad
  \rho=h_0(R)-h_0(\widetilde R)
       +h_0(\widetilde R/I)-h_0(R/I);
\]
see \cite[Exercise~III.4.4]{weibel}.
Since \(K_{-1}(X)\)=0, we must have \(\rho=0\). Hence
\[
h_0(R)-h_0(\widetilde{R})
+h_0(\widetilde{R}/I\widetilde{R})
-h_0(R/I)=0.
\]
Therefore,
\begin{equation}\label{h0R}
h_0(\widetilde{R}/I\widetilde{R})=h_0(R/I). 
\end{equation}

Let
\[
Z:=\operatorname{Spec}(R/I).
\]
The singular locus of \(X\) is supported on \(Z\). For the normalization morphism
\[
\nu\colon \operatorname{Spec}(\widetilde{R})\longrightarrow X,
\]
we have
\[
\nu^{-1}(Z)=
\operatorname{Spec}(\widetilde{R}/I\widetilde{R}).
\]
By (\ref{h0R}), \(Z\) and \(\nu^{-1}(Z)\) have the same number of connected components. Consequently, every singular point of \(X\) has at most one point lying above it under the normalization morphism. Thus, \(X\) has only unibranch singularities.
\end{proof}

    \begin{cor}\thlabel{A1}
        Let $k$ be a perfect field and $X$ be an $\A^1$-contractible, geometrically integral affine $k$-curve. Then the normalization of $X$ is $\A_k^1$ and $X$ is isomorphic to $\A_k^1$ as a Nisnevich sheaf of sets on $Sm_k$.
    \end{cor}
    \begin{proof}
 Since $X$ is $\A^1$-contractible, $\overline{X}=X \times_k \overline{k}$ is also $\A^1$-contractible and thus $\overline{X}$ must be $\A^1$-
connected. Now note that $\overline{X}$ is not $\A^1$-invariant as a Nisnevich sheaf of sets on $Sm_{\overline{k}}$. Thus there exists $W \in Sm_{\overline{k}}$ and a morphism $H: W \times \A_{\overline{k}}^1 \to \overline{X}$, 
such that $H(\underline{\hspace{0.35em}},0) \neq H(\underline{\hspace{0.35em}},1): W \to \overline{X}$. 
Therefore, we get a non constant morphism $\tilde{H}: \A_{\overline{k}}^1 \to \overline{X}$. 
Since $\mathrm{dim}\,X=1$, $\tilde{H}$ is a dominant morphism and hence we get a dominant morphism to the normalization $\widetilde{\overline{X}}$ of $\overline{X}$ given by $H_1: \A_{\overline{k}}^1 \to \widetilde{\overline{X}}$. Therefore, $\widetilde{\overline{X}} \cong \A_{\overline{k}}^1$. Now since $k$ is a perfect field, it follows that the normalization $\widetilde{X}$ of $X$ is $\A_k^1$.
Now since $X$ is $\A^1$-contractible, by \thref{uni-branched}, it has at most unibranched singularities and hence the normalization map $\A_k^1 \to X$ is universally bijective. Therefore, on $Sm_k$ the induced map of Nisnevich sheaves $h_{\A_k^1} \to h_X$ is an isomorphism. 

    \end{proof}

    \begin{rem}\thlabel{rem3.3}
    A classical result of Lin-Zaidenberg states that over an algebraically closed field of characteristic zero, a rational affine plane curve with unibranch singularity and one place at infinity can be represented as $z^r+t^s=0$, with $\gcd(r,s)=1$, up to an automorphism of the ambient space $\A^2$ (\cite{LZ}).
        \thref{uni-branched} shows that an integral affine $\A^1$-contractible plane curve satisfies these properties and thus can be represented in the above form. 
    \end{rem}

    \subsection{Sufficient conditions for stable $\A^1$-contractibility of Koras--Russell type varieties.}

We first prove a necessary and sufficient condition for $\A^1$-contractibility of a subfamily of the smooth varieties given by the rings $A$ of the form \eqref{a}.  

    \begin{thm}\thlabel{sufficientunstable}
        Let $k$ be an algebraically closed field of characteristic zero. Let
        $$
        A= \dfrac{k[x_1,\ldots,x_m,y,z,t]}{(x_1^{r_1}\cdots x_m^{r_m}y+f(z,t)+x_m)} \text{~with~} r_i>1,
        $$
        and assume that $\mathrm{Spec}\,A$ is smooth. Then $\mathrm{Spec}\,A$ is $\A^1$-contractible if and only if $\mathrm{Spec}\,(k[z,t]/(f))$ is an integral $\A^1$-contractible plane curve.
    \end{thm}
   \begin{proof}
       If $\mathrm{Spec}\,A$ is $\A^1$-contractible, then $A$ is a UFD and $A^*=k^*$. Therefore, by \thref{thm1}, $f(z,t)$ is irreducible in $k[z,t]$. 
       Now the necessary part follows from \thref{thm3}(b).
       Conversely, if $\mathrm{Spec}\,(k[z,t]/(f))$ is a smooth $\A^1$-contractible curve, then $k[z,t]/(f)=k^{[1]}$, and hence $k[z,t]=k[f]^{[1]}$ (cf. \cite{AM}). Therefore, $A=k^{[m+2]}$, in particular $\mathrm{Spec}\,A$ is $\A^1$-contractible.
       On the other hand, if $\{f(z,t)=0\}$ is not smooth, then by \thref{rem3.3}, $f(z,t)=z^r+t^s$, with $\mathrm{gcd}(r,s)=1$, up to some change of coordinates in $k[z,t]$. Therefore, the result follows from \cite[Theorem 1.1]{DF}, \cite[Proposition 5]{DG}. 
   \end{proof}
We now consider the larger family of smooth varieties given by rings $A$ of the form \eqref{a}, i.e.,
$$
 A:=\dfrac{k[x_1,\ldots,x_m,y,z,t]}{\left( a(x_m)b(x_1,\ldots,x_{m-1})y+f(z,t)+x_m\right)},
$$
and give a sufficient condition for $\mathrm{Spec}\,A$ to be stably $\A^1$-contractible, over algebraically closed fields of characteristic zero. The family of varieties contains the Koras-Russell threefolds given by rings 
$$
B_1=\dfrac{k[x_1,y,z,t]}{\left( x^{r}y+z^s+t^u+x_1 \right)}, \text{~with~} r,s,u>1, \mathrm{gcd}(s,u)=1.
$$
 Thus the result generalizes the stable $\A^1$-contractibility of Koras--Russell threefolds proved by Hoyois, Krishna, and \O stv\ae r (\cite{HKO}). 
We first recall a result from \cite{HKO} which will be important to prove the sufficient condition. The following version can be found in \cite[Theorem 3.1]{DPO}. 

\begin{thm}\thlabel{chowgp}
Suppose $X$ is a smooth affine $k$-variety over an algebraically closed field of characteristic zero. If, for every smooth affine $k$-variety $Y$, the natural map
$
X \times Y \longrightarrow Y
$
induces isomorphisms on higher Chow groups, then $X$ is stably $\mathbb{A}^1$-contractible, in the sense that it becomes $\A^1$-contractible after taking suspension by finitely many copies of $\mathbb{P}^1$.
    
\end{thm}

\begin{thm}\thlabel{sufficient}
    
    Let $k$ be an algebraically closed field of characteristic zero, $A$ be as in \eqref{a} such that $\mathrm{Spec}\,A$ is a smooth affine $k$-variety. If for every root $\lambda$ of $a(x_m)$, $\mathrm{Spec}\left( k[z,t]/(f+\lambda) \right)$ is an integral $\A^1$-contractible curve in $\mathcal{H}((Sm_k)_{Nis}, \A^1)$, then $\mathrm{Spec}\,A$ is stably $\A^1$-contractible; that is, $\mathrm{Spec} \,A$ is $\A^1$-weakly equivalent to $\mathrm{Spec} k$ in $\mathcal{H}((Sm_k)_{Nis}, \A^1)$ after finitely many suspensions by $\mathbb{P}_k^1$.

\end{thm}

\begin{proof}

    Let $C_{\lambda}=k[z,t]/(f+\lambda)$. Since $\mathrm{Spec}\, C_{\lambda}$ is $\A^1$-contractible,  by \thref{A1}, the normalization of $C_{\lambda}$ must be $\widetilde{C_{\lambda}}:=k[u]=(k^{[1]})$. Let $I_{\lambda}$ denote the conductor ideal of $k[u_{\lambda}]$ into $C_{\lambda}$. Let $I_{\lambda}=(p(u))\widetilde{C_{\lambda}}$. Since $\mathrm{Spec}\, C_{\lambda}$ has at most unibranch singularity (by \thref{uni-branched}), 
    $\mathrm{Spec}\,(\widetilde{C_{\lambda}}/(p(u)))$ and $\mathrm{Spec}\,(C_{\lambda}/p(u))$ have equal number of connected components. 
    Therefore, for any affine $k$-algebra $D$, the inclusion $C_{\lambda} \hookrightarrow \widetilde{C_{\lambda}}$ induces the following isomorphism of total higher Chow groups:
    
    \begin{equation}\label{sing}
    \mathrm{CH}_*(D \otimes_k \widetilde{C_{\lambda}}/(p(u)), i) 
    \xrightarrow{\cong} 
    \mathrm{CH}_*(D \otimes_k C_{\lambda}/(p(u)),i),
    \end{equation}
 for every $i \geqslant 0$. 

 We now consider the following commutative diagram of total higher Chow groups:

\smallskip

  \adjustbox{scale=0.70,center}{
 		\begin{tikzcd}[row sep=1.5cm, column sep=1cm]
			 \arrow{r} & \mathrm{CH}_{*}(D \otimes_k\widetilde{C_{\lambda}},i+1) \arrow{r} \arrow{d}
			& \mathrm{CH}_{*}((D \otimes_k\widetilde{C_{\lambda}})[p(u)^{-1}],i+1) \arrow{r} \arrow{d}{\cong} 
			& \mathrm{CH}_{*}((D\otimes_k\widetilde{C_{\lambda}})/(p(u)),i) \arrow{r} \arrow{d}{\cong} 
			& \mathrm{CH}_{*}((D\otimes_k\widetilde{C_{\lambda}}),i) \arrow{r} \arrow{d}
			& \mathrm{CH}_*((D\otimes_k\widetilde{C_{\lambda}})[ (p(u))^{-1}],i) \arrow{d}{\cong} \\
			 \arrow{r}& \mathrm{CH}_{*}(D \otimes_k C_{\lambda},i+1) \arrow{r}  & \mathrm{CH}_{*}((D \otimes_k C_{\lambda})[p(u)^{-1}],i+1)\arrow{r} & \mathrm{CH}_{*}((D \otimes_k C_{\lambda})/(p(u)),i) \arrow{r} &\mathrm{CH}_*(D \otimes_k C_{\lambda},i) \arrow{r} &\mathrm{CH}_*((D \otimes_k C_{\lambda})[(p(u))^{-1}],i).
		\end{tikzcd}
		}

\smallskip

        From the above diagram, using the fact that $\widetilde{C_\lambda}=k[u]$, we obtain the following isomorphisms;

    \begin{equation}\label{8}
    \mathrm{CH}_*(D \otimes_k\widetilde{C_{\lambda}},i)=\mathrm{CH}_*(D[u],i)=\mathrm{CH}_*(D,i) \xrightarrow{\cong} \mathrm{CH}_*(D \otimes_k C_{\lambda},i), \text{~for every}\, i \geqslant 0.
\end{equation}

      We now consider the ring $A$ and the flat morphism $S:=k[x_1,\ldots,x_m] \hookrightarrow A$ (cf. \cite[Corollary 2, 20.F]{Matc}). Let $q_1,\ldots,q_l$ and $p_1,\ldots,p_n$ be the prime factors of $a(x_m)$ and $b(x_1,\ldots,x_{m-1})$ respectively.
      Now for a $k$-algebra $B$, let $S_B:=B[x_1,\ldots,x_m]$, $A_B:=A \otimes_k B$, $R=B[x_1,\ldots,x_{m-1}]$. Note that $S_B=R[x_m] \hookrightarrow A_B$ is a flat morphism. We now prove that this morphism induces isomorphisms on higher Chow groups. We first note that the morphism $S_B \hookrightarrow A_B$ induces the morphism

      \begin{align}\label{repeat}
      &S_B[p_1^{-1},\ldots,p_{i-1}^{-1}]/(p_i) =R[p_1^{-1},\ldots,p_{i-1}^{-1}][x_m]/(p_i) \to A_B[p_1^{-1},\ldots,p_{i-1}^{-1}]/(p_i) =R[p_1^{-1},\ldots,p_{i-1}^{-1}][y,z,t]/(p_i), \nonumber \\
 &\text{which induces isomorphisms on higher Chow groups.} 
 \end{align}

         We now consider the following commutative diagram of localization sequences of higher Chow groups with respect to the element $p_1$:
       \[
		\begin{tikzcd}[row sep=large, column sep=tiny]
			 \arrow{r} & \mathrm{CH}_{*}(S_B,i) \arrow{r} \arrow{d} 
			& \mathrm{CH}_{*}(S_B[p_1^{-1}],i) \arrow{r} \arrow{d}
			& \mathrm{CH}_{*}\left(\dfrac{S_B}{(p_1)},i-1\right) \arrow{r} \arrow{d} {\cong}
			& \mathrm{CH}_{*}(S_B,i-1) \arrow{r} \arrow{d} 
			& \mathrm{CH}_*(S_B[p_1^{-1}],i-1) \arrow{d} \\
			 \arrow{r}& \mathrm{CH}_{*}(A_B,i) \arrow{r}  & \mathrm{CH}_{*}(A_B[p_1^{-1}],i)\arrow{r} & \mathrm{CH}_{*}\left(\dfrac{A_B}{(p_1)},i-1\right) \arrow{r} &\mathrm{CH}_*(A_B,i-1) \arrow{r} &\mathrm{CH}_*(A_B[p_1^{-1}],i-1).
		\end{tikzcd}
		\]

        From the above diagram it is clear that for every $i \geqslant 0$, $\mathrm{CH}_*(S_B,i) \cong \mathrm{CH}_*(A_B,i)$ if and only if $\mathrm{CH}_*(S_B[p_1^{-1}],i) \cong \mathrm{CH}_*(A_B[p_1^{-1}],i)$. In the next step considering the commutative diagram of localization sequences induced by the flat morphism $S_B[p_1^{-1}] \hookrightarrow A_B[p_1^{-1}]$ with respect to the element $p_2$, we get that $\mathrm{CH}_*(S_B[p_1^{-1}],i) \cong \mathrm{CH}_*(A_B[p_1^{-1}],i)$ if and only if $\mathrm{CH}_*(S_B[p_1^{-1},p_2^{-1}],i) \cong \mathrm{CH}_*(A_B[p_1^{-1},p_2^{-1}],i)$, using the isomorphisms in \eqref{repeat}. Repeating this process for all the prime factors of $b(x_1,\ldots,x_{m-1})$, we get that $\mathrm{CH}_*(S_B,i) \cong \mathrm{CH}_*(A_B,i)$ if and only if $\mathrm{CH}_*(S_B[p_1^{-1},\ldots,p_n^{-1}],i) \cong \mathrm{CH}_*(A_B[p_1^{-1},\ldots,p_n^{-1}],i)$. 

        Now let $q=(q_1\cdots q_l) \in k[x_m]$. Considering the following commutative diagram of localization sequences induced by the flat morphism $S_B[p_1^{-1},\ldots,p_n^{-1}] \hookrightarrow A_B[p_1^{-1},\ldots,p_n^{-1}]$ with respect to the element $q$, we get

        \smallskip
\adjustbox{scale=0.70,center}{
         \begin{tikzcd}[row sep=1.2cm, column sep=0.6cm]
			 \arrow{r} & \mathrm{CH}_{*}(S_B[p_1^{-1},\ldots,p_n^{-1}],i) \arrow{r} \arrow{d} 
			& \mathrm{CH}_{*}(S_B[p_1^{-1},\ldots,p_n^{-1},q^{-1}],i) \arrow{r} \arrow{d}{\cong}
			& \mathrm{CH}_{*}\left(\dfrac{S_B[p_1^{-1},\ldots,p_n^{-1}]}{(q)},i-1\right) \arrow{r} \arrow{d} 
			& \mathrm{CH}_{*}(S_B[p_1^{-1},\ldots,p_n^{-1}],i-1) \arrow{r} \arrow{d}
			& \mathrm{CH}_*(S_B[p_1^{-1},\ldots,p_n^{-1},q^{-1}],i-1) \arrow{d}{\cong} \\
			 \arrow{r}& \mathrm{CH}_{*}(A_B[p_1^{-1},\ldots,p_n^{-1}],i) \arrow{r}  & \mathrm{CH}_{*}(A_B[p_1^{-1},\ldots,p_n^{-1},q^{-1}],i)\arrow{r} & \mathrm{CH}_{*}\left(\dfrac{A_B[p_1^{-1},\ldots,p_n^{-1}]}{(q)},i-1\right) \arrow{r} &\mathrm{CH}_*(A_B[p_1^{-1},\ldots,p_n^{-1}],i-1) \arrow{r} &\mathrm{CH}_*(A_B[p_1^{-1},\ldots,p_n^{-1},q^{-1}],i-1).
		\end{tikzcd}
		}

        \smallskip
        From the above diagram we have $\mathrm{CH}_*(S_B[p_1^{-1},\ldots,p_n^{-1}],i) \cong \mathrm{CH}_*(A_B[p_1^{-1},\ldots,p_n^{-1}],i)$ for every $i \geqslant 0$ if and only if $\mathrm{CH}_*\left(\dfrac{S_B[p_1^{-1},\ldots,p_n^{-1}]}{(q)},i\right) \cong \mathrm{CH}_{*}\left(\dfrac{A_B[p_1^{-1},\ldots,p_n^{-1}]}{(q)},i\right)$.
        Now let $B_1:=B[x_1,\ldots,x_{m-1},p_1^{-1},\ldots,p_n^{-1}]$. Let $\lambda_1,\ldots,\lambda_r$ be roots of $a(x_m)$ in $k$. 
        Now, from the above commutative diagram, we get the morphism
        $$
        \dfrac{S_B[p_1^{-1},\ldots,p_n^{-1}]}{(q)} = B_1[x_m]/(q) =\bigoplus_{\lambda_i} B_1 \hookrightarrow \dfrac{A_B[p_1^{-1},\ldots,p_n^{-1}]}{(q)} = \bigoplus_{\lambda_i} \dfrac{B_1[y,z,t]}{(f(z,t)+\lambda_i)} = \bigoplus_{\lambda_i} (B_1 \otimes_k C_{\lambda_i})[y], 
        $$
      induces isomorphisms on total higher Chow groups $\mathrm{CH}_*(-,i)$ for every $i \geqslant 0$, using the isomorphism in \eqref{8} for $D=B_1$.
        Therefore, $\mathrm{CH}_*(S_B[p_1^{-1},\ldots,p_n^{-1}],i) \cong \mathrm{CH}_*(A_B[p_1^{-1},\ldots,p_n^{-1}],i)$ for every $i \geqslant 0$. This implies that $\mathrm{CH}_*(S_B,i) =\mathrm{CH}_*(B,i)\cong \mathrm{CH}_*(A_B,i)$ for every $i \geqslant 0$. Therefore, for every affine $k$-variety $B$, $\mathrm{Spec}\, A \times \mathrm{Spec}\, B \to \mathrm{Spec}\,B$ induces isomorphisms $\mathrm{CH}^j(B,i) \to \mathrm{CH}^j(A_B,i)$ for every $i,j \geqslant 0$. 
Now using \thref{chowgp}, we get that $\mathrm{Spec}\, A$ is stably $\A^1$-contractible.
\end{proof}

 \section{Application to the Embedding Problem and characterization of affine spaces.}
 
 In this section we prove Theorems 1.4 and 1.5, as mentioned in the introduction and answer Question 1(ii) for a large subfamily of varieties given by $\mathrm{Spec}\,A$, as in \eqref{a}.

    \begin{thm}\thlabel{embedding1}
        Let $k$ be a field of characteristic zero, and 
        $$
        A=\dfrac{k[x_1,\ldots,x_m,y,z,t]}{(H)},
        $$
        where
        $$
        H=a(x_m)b(x_1,\ldots,x_{m-1})y+f(z,t)+x_m
        $$ such that
        $a(x_m)$ has at least two distinct roots in $\overline{k}$. Then the following are equivalent:
        \begin{enumerate}[\rm (a)]
            \item $\mathrm{Spec}\,A$ is smooth and $\A^1$-contractible.

            \item $A=k^{[m+2]}$.

            \item $k[z,t]=k[f]^{[1]}$.

            \item $k[x_1,\ldots,x_m,y,z,t]=k[H]^{[m+2]}$.
        \end{enumerate}
    \end{thm}

    \begin{proof}
        Note that $\rm (c) \Rightarrow (d) \Rightarrow (b) \Rightarrow (a)$ hold. We show $\rm (a) \Rightarrow (c)$. 
        
        Let $\overline{A}=A \otimes_k \overline{k}$. Since $\mathrm{Spec}\, A$ is smooth and $\A^1$-contractible, so is $\mathrm{Spec}\, \overline{A}$. Therefore, by \thref{thm3}, we get that for every root $\lambda$ of $a(x_m)$ in $\overline{k}$, $\overline{k}[z,t]/(f+\lambda)$ is $\A^1$-contractible. 
Let $\lambda_1,\lambda_2$ be two distinct roots of $a(x_m)$. Suppose, for contradiction, $\overline{k}[z,t]/(f+\lambda_1)$ is a singular $\A^1$-contractible plane curve. In that case, by \thref{rem3.3}, we get that $f(z,t)-\lambda_1=z^r+t^s$, for some $r,s\geqslant 2$ with $\mathrm{gcd}(r,s)=1$, up to a change of coordinates in $\overline{k}[z,t]$. But this contradicts that $\overline{k}[z,t]/(f+\lambda_2)$ is $\A^1$-contractible. Therefore, $\overline{k}[z,t]/(f+\lambda_1)$ must be smooth $\A^1$-contractible, and hence $\overline{k}[z,t]/(f+\lambda_1)=\overline{k}^{[1]}$. Since $\operatorname{char} k=0$, it follows that $\overline{k}[z,t]=\overline{k}[f]^{[1]}$ (cf. \cite{AM}), and since $k$ is a perfect field, $k[z,t]=k[f]^{[1]}$, because separable $\A^1$-forms over a PID are polynomial rings (cf. \cite{D}). 
        
    \end{proof}
    
    \begin{thm}\thlabel{embedding}
		Let $k$ be any field and $A$ be as in \eqref{b}, i.e.,
        $$
        A=\frac{k[x_1,\ldots,x_m,y,z,t]}{(H)}, 
       $$        
where 
        $$
        H = x_1^{r_1}a_1(x_1) \cdots x_{m-1}^{r_{m-1}}a_{m-1}(x_{m-1})x_m^{r_m}a_m(x_m) y + f(z, t) + x_m \text{~with~} r_i>1.
        $$ 
Then the following are equivalent:
		\begin{enumerate}[\rm (a)]
		\item $A=k^{[m+2]}$ 
		\item $k[Z,T]=k[f]^{[1]}$
		\item $k[X_1,\ldots,X_m,Y,Z,T]=k[H]^{[m+2]}$.
		\end{enumerate}
	\end{thm}
	\begin{proof}
		Note that (b) $\Rightarrow$ (c) $\Rightarrow$ (a) hold. We now show (a) $\Rightarrow$ (b). 

        Since $A \otimes_k \overline{k}=\overline{k}^{[m+2]}$, without loss of generality we can assume $f(z,t)=b_0(z)+b_1(z)t$ in $\overline{k}[z,t]$ (by \thref{thm1}), hence $\mathrm{Spec}(k[z,t]/(f))$ is a smooth curve. 
        Therefore by \thref{thm3}, we have $\overline{k}[z,t]/(f)=\overline{k}^{[1]}$ and since $f(z,t)=b_0(z)+b_1(z)t$ in $\overline{k}[Z,T]$, it follows that $\overline{k}[z,t]=\overline{k}[f]^{[1]}$.
        
        If $k$ is a perfect field, then $k[z,t]=k[f]^{[1]}$ (cf. \cite{D}).

		If $k$ is not perfect, then $k$ must be infinite and hence from $A=k^{[m+2]}$ without loss of generality we get $f(z,t)=b_0(z)+b_1(z)t$ for some $b_0(z),b_1(z) \in k[z]$ (cf. \thref{thm1}). Therefore $\overline{k}[z,t]/(f)=\overline{k}^{[1]}$ implies $k[z,t]=k[f]^{[1]}$.
	\end{proof}

    \section{A new family of smooth \(\mathbb A^1\)-contractible varieties}

\begin{thm}\label{0.2}
Let \(k\) be a  field, let \(m\ge2\), and let \(r,s>1\) be relatively
prime integers. Let
\begin{equation}\label{b.1}
A=\frac{k[x_1,\ldots,x_m,y,z,t]}{(b(x_1,\ldots,x_{m-1}) x_m^{r_m}y+z^r+t^s+x_m)},
\end{equation}
where $r_m>1$. 
Then $\operatorname{Spec}A\text{ is smooth and }\mathbb A^1\text{-contractible}$.
\end{thm}

\begin{proof}
Let
\[
R=k[x_1,\ldots,x_{m-1}],\qquad
W=\operatorname{Spec}R=\mathbb A^{m-1}_k,
\]
and consider the projection
\[
\pi:X:=\operatorname{Spec}A\longrightarrow W.
\]
We first verify that $\pi$ is smooth. Let
\[
F=b(x_1,\ldots,x_{m-1})x_m^{r_m}y+z^r+t^s+x_m.
\]
Relative to $R$, two of the partial derivatives are
\[
\frac{\partial F}{\partial y}=b\,x_m^{r_m},
\qquad
\frac{\partial F}{\partial x_m}
=r_m b\,x_m^{r_m-1}y+1.
\]
If $J$ denotes the relative Jacobian ideal, then
\[
x_m\frac{\partial F}{\partial x_m}
-r_my\frac{\partial F}{\partial y}=x_m,
\]
so $x_m\in J$. Since $r_m>1$, we also have
$r_m b\,x_m^{r_m-1}y\in (x_m)\subseteq J$, and therefore
\[
1=\frac{\partial F}{\partial x_m}
-r_m b\,x_m^{r_m-1}y\in J.
\]
Thus the relative Jacobian ideal is the unit ideal, and the Jacobian criterion
implies that $\pi$ is smooth.
Since $W$ is smooth over $k$, it follows that $X$ is smooth over $k$.

Let $\mathfrak p\in\operatorname{Spec}R$, let
$\kappa(\mathfrak p)=R_{\mathfrak p}/\mathfrak pR_{\mathfrak p}$, and let
$\overline b$ be the image of $b$ in $\kappa(\mathfrak p)$. The fiber of
$\pi$ over $\mathfrak p$ is
\[
X_{\mathfrak p}
=\operatorname{Spec}
\frac{\kappa(\mathfrak p)[x_m,y,z,t]}
{\bigl(\overline b\,x_m^{r_m}y+z^r+t^s+x_m\bigr)}.
\]
There are two cases. If $b\in\mathfrak p$, then
\[
X_{\mathfrak p}
\cong\operatorname{Spec}\kappa(\mathfrak p)[y,z,t]
=\mathbb A^3_{\kappa(\mathfrak p)},
\]
which is $\mathbb A^1$-contractible. If $b\notin\mathfrak p$, then
$\overline b\in\kappa(\mathfrak p)^*$ and, after the change of variable
$Y=\overline b\,y$,
\[
X_{\mathfrak p}\cong
\operatorname{Spec}
\frac{\kappa(\mathfrak p)[x_m,Y,z,t]}
{\bigl(x_m^{r_m}Y+z^r+t^s+x_m\bigr)}.
\]
This is $\mathbb A^1$-contractible over
$\kappa(\mathfrak p)$ by \cite{DF} and 
\cite[Proposition~4.15, Remark~4.18]{Vij}.

Now since $\pi$ is smooth and that every residue-field
fiber $X_{\mathfrak p}$ is $\mathbb A^1$-contractible, therefore by,
\cite[Proposition~1.15]{DMO2026} 
\[
\pi:X\longrightarrow W
\]
is an $\mathbb A^1$-weak equivalence. Since
$W=\mathbb A^{m-1}_k$ is $\mathbb A^1$-contractible, so is $X$.
\end{proof}

We now consider the following subfamily of varieties given by \eqref{b.1} 
and show that they are not isomorphic to the $\mathbb{A}^1$-contractible varieties in
 \cite{DG} and \cite{Vij} given by
\begin{equation}\label{c1}
      \mathcal{X}(\underline{r},p):=\left(\frac{k[x_1,\ldots,x_m,y,z,t]}
 {\bigl(x_1^{r_1}x_2^{r_2}\cdots x_m^{r_m}y+z^r+t^s+x_m p(x_1\cdots x_m)\bigr)}\right),  
\end{equation}
where $r_i>1$ and $p(0) \neq 0$.

Let $a(x_1) \in k[x_1] \setminus k$ satisfy \(a(0)\ne0\), and assume that it has only multiple roots in an algebraic closure $\overline{k}$ of $k$. For \(m\ge2\) and relatively
prime integers \(r,s>1\), consider the variety
\begin{equation}\label{c}
\mathcal{X}_a:=\operatorname{Spec}
\frac{k[x_1,\ldots,x_m,y,z,t]}
 {\bigl(x_1^{r_1}a(x_1)x_2^{r_2}\cdots x_m^{r_m}y+z^r+t^s+x_m\bigr)},
\end{equation}
where $r_i>1$ for all $i$.
By Theorem~\ref{0.2}, $\mathcal{X}_a$ is smooth and \(\mathbb A^1\)-contractible. 
We now show that $\mathcal{X}_a$ is not isomorphic to any of the varieties $\mathcal{X}(\underline{r},p)$ or affine spaces, therefore they provide new examples of \(\mathbb A^1\)-contractible varieties which are not affine spaces.
The next results compute the Derksen and Makar--Limanov invariants of \(\mathcal{X}_a\) over an algebraically closed field.

\begin{prop}\label{0.3}
        Let \(k\) be an algebraically closed field, and $\mathcal{X}_a$ be as in \eqref{c}. 
        Then the Derksen invariant $\mathrm{DK}(\mathcal{X}_a)=k[x_1,\ldots,x_m,z,t]$.
\end{prop}
\begin{proof}
 Note that \(\mathcal{X}_a= \mathrm{Spec}\,A_a\) where
$$
 A_a:=\frac{k[x_1,\ldots,x_m,y,z,t]}
 {\bigl(x_1^{r_1}a(x_1)x_2^{r_2}\cdots x_m^{r_m}y+z^r+t^s+x_m\bigr)}.
$$ 
 Since $k$ is algebraically closed,
\[
a(x_1)=c\prod_{i=1}^n(x_1-\lambda_i)^{e_{\lambda_i}},
\qquad c\in k^*,\quad e_{\lambda_i} >1.
\]
Let $R_i=k[x_1,\ldots,x_{i-1},x_{i+1},\ldots,x_m]$ for $2 \leqslant i \leqslant m$.
Note that
\begin{equation}\label{e}
A_a=\frac{R_i[x_i,y,z,t]}{(s_i x_i^{r_i}y+z^r+t^s+x_m)},
\end{equation}
where $s_i=x_1^{r_1}a(x_1)x_2^{r_2}\cdots x_{i-1}^{r_{i-1}}x_{i+1}^{r_{i+1}}\cdots x_m^{r_m}$, for $i \neq 1$.
Further
\begin{equation}\label{f}
A_a=\frac{R_1[x_1,y,z,t]}{(s_1 x_1^{r_1}y+z^r+t^s+x_m)},
\end{equation}
where $s_1=a(x_1)x_2^{r_2}\cdots x_m^{r_m}$ and
\begin{equation}\label{g}
A_a=\frac{R_1[x_1,y,z,t]}{(s_{\lambda_j}(x_1-\lambda_{j})^{e_{\lambda_j}}y+z^r+t^s+x_m)},
\end{equation}
where $s_{\lambda_j}=cx_1^{r_1}\prod_{i\neq j}(x_1-\lambda_i)^{e_{\lambda_i}}x_2^{r_2}\cdots x_m^{r_m}$, for $1 \leqslant j \leqslant n$.
For \(1\le i\le m\), consider the degree function on $A_a$
\[ 
\omega_i(x_i)=-1,\qquad \omega_i(y)=r_i,\qquad
\omega_i(x_j)=\omega_i(z)=\omega_i(t)=0\quad(j\ne i).
\]
By \cite[Theorem~5.5]{pre} the degree functions induce admissible $\mathbb{Z}$-filtrations on \(A_a\) (\cite[Definition~5.3]{pre})
such that the associated graded rings are given as follows:
\begin{equation}\label{h.1}
\operatorname{gr}_{\omega_1}(A_a)
 \cong\frac{k[x_1,\ldots,x_m,y,z,t]}{(x_1^{r_1}x_2^{r_2}\cdots x_m^{r_m}y+z^r+t^s+x_m)}
\end{equation}
and 
\begin{equation}\label{i}
\operatorname{gr}_{\omega_i}(A_a)
 \cong\frac{k[x_1,\ldots,x_m,y,z,t]}{(x_1^{r_1}a(x_1)x_2^{r_2}\cdots x_m^{r_m}y+z^r+t^s+x_m)}
\end{equation}
for $2 \leq i < m$, and
\begin{equation}\label{j}
\operatorname{gr}_{\omega_m}(A_a)
 \cong\frac{k[x_1,\ldots,x_m,y,z,t]}{(x_1^{r_1}a(x_1)x_2^{r_2}\cdots x_m^{r_m}y+z^r+t^s)}
\end{equation}
for $i = m$.
Now for $1 \leqslant j \leqslant n$, consider the degree function on $A_a$
$$ 
{\omega}_{\lambda_j}(x_1-\lambda_j)=-1,\qquad {\omega}_{\lambda_j}(y)=e_{\lambda_j},\qquad
{\omega}_{\lambda_j}(x_i)=\omega_{\lambda_j}(z)=\omega_{\lambda_j}(t)=0\quad( i \neq 1).
$$
Again by \cite[Theorem~5.5]{pre}, these degree functions induce admissible $\mathbb{Z}$-filtrations 
on \(A_a\) such that the associated graded rings are given as follows:
\begin{equation}\label{k.1}
\operatorname{gr}_{\omega_{\lambda_j}}(A_a)
 \cong\frac{k[x_1,x_2,\ldots,x_m,y,z,t]}
        {\bigl((x_1-\lambda_j)^{e_{\lambda_j}}
          \prod_{i=2}^m x_i^{r_i}y+z^r+t^s+x_m\bigr)},
\end{equation}
Let \(w\) be any one of the degree functions above. Suppose that a
nontrivial exponential map \(\varphi\) on \(A_a\) fixes a nonzero
\(g\) with \(w(g)>0\). Let $\rho_{w}: A_a \to \mathrm{gr}_w(A_a)$ be the natural projection map.
Then $\varphi$ induces a nontrivial exponential map \(\overline\varphi\) on \(\operatorname{gr}_w(A_a)\) such that
\[
\rho_w(g)\in
 \operatorname{gr}_w(A_a)^{\overline\varphi},
\qquad
\rho_w(y)\mid \rho_w(g),
\]
where \(\overline\varphi\) is nontrivial
(\cite[Theorem~5.4 and Remark~5.12(iii)]{pre}).
Therefore, $\rho_w(y)$ is invariant under \(\overline\varphi\).
Hence there exist nontrivial exponential maps on the rings as in 
\eqref{h.1}, \eqref{i}, \eqref{j} and \eqref{k.1} which fix the images of \(y\) in the respective rings.
But this contradicts \cite[Theorem~5.23]{pre}(ii).
Therefore, for every element $g$ in $\mathrm{DK}(\mathcal{X}_a)$, the degree of $g$ with respect to each degree function $\omega_i$, $1 \leqslant i \leqslant m$,
and $\omega_{\lambda_j}$, for $1 \leqslant j \leqslant n$, is nonpositive. 

Therefore,
\begin{align*}
        \mathrm{DK}(\mathcal{X}_a) &\subseteq \bigcap_{i=1}^m k[x_1,\ldots,x_m, (x_i)^{r_i}y,z,t] \bigcap_{j=1}^n k[x_1,\ldots,x_m,(x_1-\lambda_j)^{e_{\lambda_j}}y,z,t]\\
        &=  \bigcap_{i=1}^m k[x_1,\ldots,x_m, \dfrac{z^r+t^s+x_m}{s_i},z,t] \bigcap_{j=1}^n k[x_1,\ldots,x_m,\frac{z^r+t^s+x_m}{s_{\lambda_j}},z,t]\\
        &=k[x_1,\ldots,x_m,z,t].
\end{align*}
Moreover, $k[x_1,\ldots,x_m,z,t]\subseteq \mathrm{DK}(\mathcal{X}_a)$  (\cite[Lemma~3.1(i)]{GPproc}).
This proves the result.
\end{proof}

\begin{prop}\label{0.4}
 Let \(k\) be an algebraically closed field, and $\mathcal{X}_a$ be as in \eqref{c}. 
        Then the Makar-Limanov invariant $\mathrm{ML}(\mathcal{X}_a)=k[x_1,\ldots,x_m]$.
\end{prop}
\begin{proof}
Let $A_a$ be the ring as in Proposition~\ref{0.3}. By \cite[Lemma~3.1(ii)]{GPproc},
\[
\operatorname{ML}(\mathcal{X}_a)\subseteq k[x_1,\ldots,x_m].
\]
Suppose, for contradiction, \(x_i \notin \mathrm{ML}(\mathcal{X}_a)\). Then there exists a nontrivial
exponential map \(\varphi\) such that \(x_i \notin A_a^\varphi\). Since
\(\operatorname{tr.deg}_k A_a^\varphi=m+1\), we choose algebraically independent elements
\[
f_1,\ldots,f_{m+1}\in A_a^\varphi \subset \mathrm{DK}(\mathcal{X}_a)=k[x_1,\ldots,x_m,z,t].
\]
Let \(E_i=k[x_1,\ldots,x_{i-1},x_{i+1},\ldots,x_m,z,t]\) and 
\[
f_j=g_j+x_i h_j,\qquad g_j\in E_i,\quad h_j\in k[x_1,\ldots,x_m,z,t].
\]
The elements \(g_1,\ldots,g_{m+1}\) are algebraically independent,
otherwise there is a nonzero polynomial \(H\) with \(H(g_1,\ldots,g_{m+1})=0\)
and that would give
\[
0\ne H(f_1,\ldots,f_{m+1})\in x_ik[x_1,\ldots,x_m,z,t]\cap A_a^\varphi.
\]
Factorial closedness of $A_a^{\varphi}$ would imply \(x_i\in A_a^\varphi\), which contradicts our assumption.
Now with respect to degree function \(\omega_i\) on \(A_a\) given by
$$
\omega_i(x_i)=-1,\qquad \omega_i(y)=r_i,\qquad\omega_i(z)=\omega_i(t)=\omega_i(x_j)=0,(~j\ne i~),
$$ 
$\varphi$ induces a nontrivial exponential map \(\widehat\varphi\) on $\mathrm{gr}_{\omega_i}(A_a)$
such that 
\begin{equation}\label{l}
        \mathrm{gr}_{\omega_i}(A_a) \cong \frac{k[x_1,\ldots,x_m,y,z,t]}{(x_1^{r_1}a(x_1)x_2^{r_2}\cdots x_m^{r_m}y+z^r+t^s+x_m)} \text{~if~} i \neq 1,m, 
\end{equation}
\begin{equation}\label{m}
        \mathrm{gr}_{\omega_i}(A_a) \cong \frac{k[x_1,\ldots,x_m,y,z,t]}{(x_1^{r_1}x_2^{r_2}\cdots x_m^{r_m}y+z^r+t^s+x_m)} \text{~if~} i = 1,
\end{equation}
\begin{equation}\label{n}
        \mathrm{gr}_{\omega_i}(A_a) \cong \frac{k[x_1,\ldots,x_m,y,z,t]}{(x_1^{r_1}a(x_1)x_2^{r_2}\cdots x_m^{r_m}y+z^r+t^s)} \text{~if~} i = m,
\end{equation}
and $\rho_i(f_j) \in \mathrm{gr}_{\omega_i}(A_a)^{\hat{\varphi}}$, where $\rho_i: A_a \to \mathrm{gr}_{\omega_i}(A_a)$ is the natural projection map.
Now note that $\rho_i(f_j)=\overline{g_j}$ where $\overline{g_j}$ is the image of $g_j$ in $\mathrm{gr}_{\omega_i}(A_a)$.
Since $\mathrm{tr.\,deg}_k k[\overline{g_1},\ldots,\overline{g_{m+1}}] = m+1$, $k[\overline{g_1},\ldots,\overline{g_{m+1}}] \subseteq k[x_1,\ldots,x_{i-1},x_{i+1},\ldots,x_m,z,t]$ has the same transcendence degree as the ambient polynomial ring. Since the ring $\mathrm{gr}_{\omega_i}(A_a)^{\hat{\varphi}}$ is algebraically closed, we get that $k[x_1,\ldots,x_{i-1},x_{i+1},\ldots,x_m,z,t] \subseteq \mathrm{gr}_{\omega_i}(A_a)^{\hat{\varphi}}$.
Therefore, from \eqref{l}, \eqref{m} and \eqref{n}, we have $y, x_i \in \mathrm{gr}_{\omega_i}(A_a)^{\hat{\varphi}}$ as $\mathrm{gr}_{\omega_i}(A_a)^{\hat{\varphi}}$ is factorially closed.
But this contradicts the nontriviality of \(\hat{\varphi}\).
Thus the initial assumption is false, and hence
$$\mathrm{ML}(\mathcal{X}_a)=k[x_1,\ldots,x_m].$$
\end{proof}

\begin{cor}\label{0.5}
        The varieties $\mathcal{X}_a$ as in \eqref{c} are not isomorphic to affine spaces.
\end{cor}
\begin{proof}
Suppose, for contradiction, that
$\mathcal{X}_a\cong\mathbb A_k^{m+2}$. After base change to an algebraic
closure $\overline{k}$ of $k$, we obtain
\[
\mathcal{X}_a\times_k\overline{k}
\cong\mathbb A_{\overline{k}}^{m+2}.
\]
Hence
\[
\mathrm{ML}(\mathcal{X}_a\times_k\overline{k})=\overline{k}.
\]
On the other hand, Proposition~\ref{0.4}, applied over $\overline{k}$, gives
\[
\mathrm{ML}(\mathcal{X}_a\times_k\overline{k})
=\overline{k}[x_1,\ldots,x_m],
\]
a contradiction. Therefore $\mathcal{X}_a$ is not isomorphic to affine space.
\end{proof}

\begin{prop}\label{0.6}
          The varieties $\mathcal{X}_a$ as in \eqref{c} are not isomorphic to any of the varieties $\mathcal{X}(\underline{r},p)$ as in \eqref{c1}.
\end{prop}
\begin{proof}

If possible, suppose that, some $\mathcal{X}_a$ is isomorphic to a
variety $\mathcal{X}(\underline r,p)$. After extending scalars to an algebraic
closure, the two varieties remain isomorphic, so we may assume that $k$ is
algebraically closed. Let
\[
A_a=
\frac{k[x_1,\ldots,x_m,y,z,t]}
{\bigl(x_1^{r_1}a(x_1)x_2^{r_2}\cdots x_m^{r_m}y+z^r+t^s+x_m\bigr)}
\]
and
\[
B=
\frac{k[x'_1,\ldots,x'_m,y',z',t']}
{\bigl((x'_1)^{n_1}\cdots(x'_m)^{n_m}y'
 +(z')^{r'}+(t')^{s'}+x'_mP(x'_1,\ldots,x'_m)\bigr)}.
\]
Let $\phi:A_a\to B$ be an isomorphism. By Propositions~\ref{0.3}
and~\ref{0.4},
\[
\operatorname{ML}(A_a)=k[x_1,\ldots,x_m]=:M_{A_a},
\qquad
\operatorname{DK}(A_a)=k[x_1,\ldots,x_m,z,t]=:D_{A_a},
\]
whereas \cite[Proposition~2.10]{adv2}, \cite[Proposition ~3.4(i)]{adv2} give
\[
\operatorname{ML}(B)=k[x'_1,\ldots,x'_m]=:M_B,
\qquad
\operatorname{DK}(B)=k[x'_1,\ldots,x'_m,z',t']=:D_B.
\]
Consequently $\phi$ restricts to isomorphisms
\[
\phi:M_{A_a}\xrightarrow{\sim}M_B,
\qquad
\phi:D_{A_a}\xrightarrow{\sim}D_B.
\]

Let
\[
\alpha=x_1^{r_1}a(x_1)x_2^{r_2}\cdots x_m^{r_m},
\qquad
q=z^r+t^s+x_m.
\]
Then $\alpha y=-q$ in $A_a$. Since
\[
B\subset D_B[(x'_1\cdots x'_m)^{-1}],
\]
there exists $N>0$ such that
\[
(x'_1\cdots x'_m)^N\phi(y)\in D_B.
\]
Using $\phi(y)=-\phi(q)/\phi(\alpha)$, we obtain
\[
\phi(\alpha)\mid (x'_1\cdots x'_m)^N\phi(q)
\qquad\text{in the UFD }D_B.
\]
In the UFD $D_{A_a}$, the elements $\alpha$ and $q$ are relatively prime. Hence
$\phi(\alpha)$ and $\phi(q)$ are relatively prime in $D_B$.

Because $a(0)\neq0$ and $a$ is nonconstant, every root of $a$ is nonzero.
We now choose a root $\lambda\in k^*$ of $a$. Both $x_1$ and $x_1-\lambda$ are
irreducible factors of $\alpha$. Therefore there exist
$j,\ell\in\{1,\ldots,m\}$ and $\mu,\mu_1\in k^*$ such that
\[
\phi(x_1)=\mu x'_j,
\qquad
\phi(x_1-\lambda)=\mu_1x'_\ell, \text{~with~} \ell \neq j.
\]
which is impossible. Thus we get a contradiction and hence the result follows.
\end{proof}

\begin{rem}
        Over a  field, in every dimension $\geqslant 4$, the varieties $\mathcal{X}_a$ give examples
        of smooth $\mathbb{A}^1$-contractible varieties (by Proposition~\ref{0.2}) that are not isomorphic to affine spaces (by Corollary~\ref{0.5})
        and are not isomorphic to the examples $\mathcal{X}(\underline{r},f)$ as given in \cite{DG} and \cite{Vij} (by Proposition~\ref{0.6}).
\end{rem}

 \section*{Acknowledgement}
  The author is thankful to Anand Sawant, Amit Hogadi, Utsav Choudhury, Marco Schlichting and Charanya Ravi for fruitful discussions and giving valuable suggestions.

\end{document}